\documentclass[11pt]{article}
\usepackage[utf8]{inputenc}
\usepackage[english]{babel}
\usepackage{graphicx}
\usepackage{amsmath}
\usepackage{amssymb}
\usepackage{amsfonts}
\usepackage{amsthm}
\usepackage[left=2.1cm,top=1.5cm,right=2.1cm, bottom=2.2cm,letterpaper]{geometry}
\usepackage{mathrsfs}
\usepackage{caption}
\usepackage{subcaption}
\usepackage{latexsym}
\usepackage{xfrac}
\usepackage{tcolorbox}
\usepackage{enumitem}
\usepackage{float}
\usepackage{hyperref}
\hypersetup{
	colorlinks=true,
	linkcolor=red,
	filecolor=magenta,
	urlcolor=cyan,
}

\usepackage[all]{hypcap}
\usepackage{autonum}
 \usepackage{url}
\usepackage[toc]{appendix}
\newtheorem{theorem}{Theorem}[section]

\newtheorem{lemma}{Lemma}[section]
\newtheorem{definition}{Definition}[section]

\newtheorem{remark}{Remark}[section]
\numberwithin{equation}{section}
\numberwithin{figure}{section}

\newcommand{\R}{\mathbb{R}}
\newcommand{\eps}{\varepsilon}

\newcommand{\bs}{\boldsymbol}
\makeatletter
\renewcommand*\env@matrix[1][*\c@MaxMatrixCols c]{%
	\hskip -\arraycolsep
	\let\@ifnextchar\new@ifnextchar
	\array{#1}}
\makeatother

\def\neweq#1{\begin{equation}\label{#1}}
\def\endeq{\end{equation}}

\begin{document}

\title{\vspace{-0.7cm}On the stationary Navier-Stokes equations in distorted pipes under energy-stable outflow boundary conditions}

\author{Alessio Falocchi -- Ana Leonor Silvestre -- Gianmarco Sperone}
\date{}
\maketitle
\vspace*{-6mm}

\begin{abstract}
	\noindent
The steady motion of a viscous incompressible fluid in distorted pipes, of finite length, is modeled through the Navier-Stokes equations with mixed boundary conditions: the inflow is given by an arbitrary member of the Lions-Magenes class with positive net influx, and the fluid motion is subject to a directional do-nothing boundary condition on the outlet, together with the standard no-slip assumption on the remaining walls of the domain. Existence of a weak solution to such Navier-Stokes system is proved without any restriction on the size of the data (that is, inlet velocity and external force) by means of the Leray-Schauder Principle, in which the required a priori estimate is obtained by a contradiction argument that employs Bernoulli's law for solutions of the stationary Euler equations. In turn, under a suitable smallness assumption on the data, we also prove the unique solvability of the boundary-value problem.
\end{abstract}
	
\section{Introduction and presentation of the problem} \label{introduction}

The study of the laminar motion of a viscous incompressible fluid in a junction of pipes $\Omega \subset \mathbb{R}^{N}$ ($N=2$ or $N=3$), whose bounding walls are rigid and impermeable, represents an essential and ubiquitous subject-matter in Fluid Mechanics \cite{landau} due to its applications in a variety of fields such as aerodynamics \cite{von2004aerodynamics}, hemodynamics \cite{galdi2008hemodynamical} and petroleum engineering \cite{bradley1987petroleum}. Mathematically \cite{ladyzhenskaya1969mathematical}, this motion is traditionally analyzed through the stationary Navier-Stokes equations in $\Omega$; they consist in determining a velocity field $\bs{u} : \Omega \longrightarrow \mathbb{R}^{N}$ and a scalar pressure $p : \Omega \longrightarrow \mathbb{R}$ satisfying the following system of Partial Differential Equations in $\Omega$:
\begin{equation} 
- \nu \Delta \bs{u} + (\bs{u} \cdot \nabla)\bs{u} + \nabla p = \bs{f} \, , \qquad \nabla \cdot \bs{u} = 0 \qquad \text{in} \ \ \ \Omega \, ,
\label{nstd0}
\end{equation}
where $\nu > 0$ denotes the kinematic viscosity of the fluid and $\bs{f} : \Omega \longrightarrow \mathbb{R}^{N}$ represents a density of body forces per unit mass. System \eqref{nstd0} is a canonical way of writing the Navier-Stokes equations when the fluid density $\varrho$ is constant, see Remark \ref{adim} for the corresponding dimensionless form. It results from scaling the Navier-Stokes equations by $\varrho$, so that the pressure is indeed the kinematic pressure \cite[Chapter 1, \S 5]{girault2012finite}. For simplicity, $p$ and $\nu$ will hereafter be referred to as the pressure and the viscosity.  
 
Now, whether the junction of pipes $\Omega$ is assumed to be unbounded or not, depends on the particular configuration intended to be modeled. On the one hand, problems arising from real-world applications rarely allow for a characterization of the fluid motion at large distances, and their numerical implementation inherently requires a bounded computational domain. On the other hand, some theoretical approaches suggest the approximation of problems, originally formulated in unbounded regions, through the introduction of consecutive \textit{invading domains}; we recall the celebrated method originally proposed by Leray in \cite{leray1933etude} for exterior problems in hydrodynamics, and successfully applied by Amick \cite{amick1977steady} and Ladyzhenskaya \& Solonnikov \cite{ladyzhenskaia1979determination} in pipes (see also the recent contribution \cite{gazzola2025steady}). In any case, a truncation of the fluid domain necessarily induces \textit{artificial boundaries} in the system of pipes, on which appropriate boundary conditions ought to be prescribed, see the articles by Blazy, Nazarov \& Specovius-Neugebauer \cite{blazy2007artificial, nazarov2008artificial} and references therein. The choice of boundary conditions to be imposed on this artificial outlet becomes a delicate question \cite{heywood1996artificial} (in both mathematical and physical terms), because the inflow is usually given, and the effects of viscosity dictate that the fluid velocity must equal zero on the walls of the system of pipes (no-slip assumption). Moreover, such boundary conditions should ensure the well-posedness of the mathematical model and the numerical stability of computational simulations; we refer to \cite{braack2014directional,bruneau2000boundary,bruneau1996new,fursikov2009optimal,kravcmar2018modeling,neustupa2022maximum,neustupa2023existence,nogueira2025regularized,nogueira2025steady,sperone2021steady} for a thorough discussion on this topic.

The families of bounded domains that will be considered in the present article are readily introduced in the next definitions:

\begin{definition} \label{addomain2}
	An open bounded set $\Omega \subset \mathbb{R}^{2}$ will be called \textbf{admissible} if $\partial \Omega$ is piecewise of class $\mathcal{C}^{2}$, $\Omega$ is simply connected and $\Omega$ is the union of three disjoint subsets as follows:
	\begin{itemize}
\item [(1)] In some coordinate system, $\Omega_{1} \doteq (0,\ell_{1}) \times (-h_{1},h_{1})$, for some $\ell_{1}>0$ and $h_{1} > 0$;
\item [(2)] In some coordinate system, $\Omega_{2} \doteq (0,\ell_{2}) \times (-h_{2},h_{2})$, for some $\ell_{2}>0$ and $h_{2} > 0$;
\item [(3)] $\Omega_{0} \doteq \Omega \setminus (\Omega_{1} \cup \Omega_{2})$  (note that $\Omega_{0}$ is not open).
	\end{itemize}
The boundary of $\Omega$ is decomposed as $ \partial \Omega = \Gamma_{I} \cup \Gamma_{W} \cup \Gamma_{O}$, where
\begin{equation}\label{boundaryomega1}
\begin{aligned}	
& \Gamma_{I} \doteq \{ 0 \} \times (-h_{1},h_{1})  \quad \text{(in the coordinate system defining $\Omega_{1}$)} \, \\[5pt]
& \Gamma_{O} \doteq \{ \ell_{2} \} \times (-h_{2},h_{2}) \quad \text{(in the coordinate system defining $\Omega_{2}$)} \, , 
\end{aligned}
\end{equation}	
and $\Gamma_{W} \subset \mathbb{R}^{2}$ represents the union of the two $\mathcal{C}^{2}$-curves connecting $\Gamma_{I}$ with $\Gamma_{O}$.
\end{definition}
\begin{figure}[H]
	\begin{center}
		\includegraphics[scale=0.65]{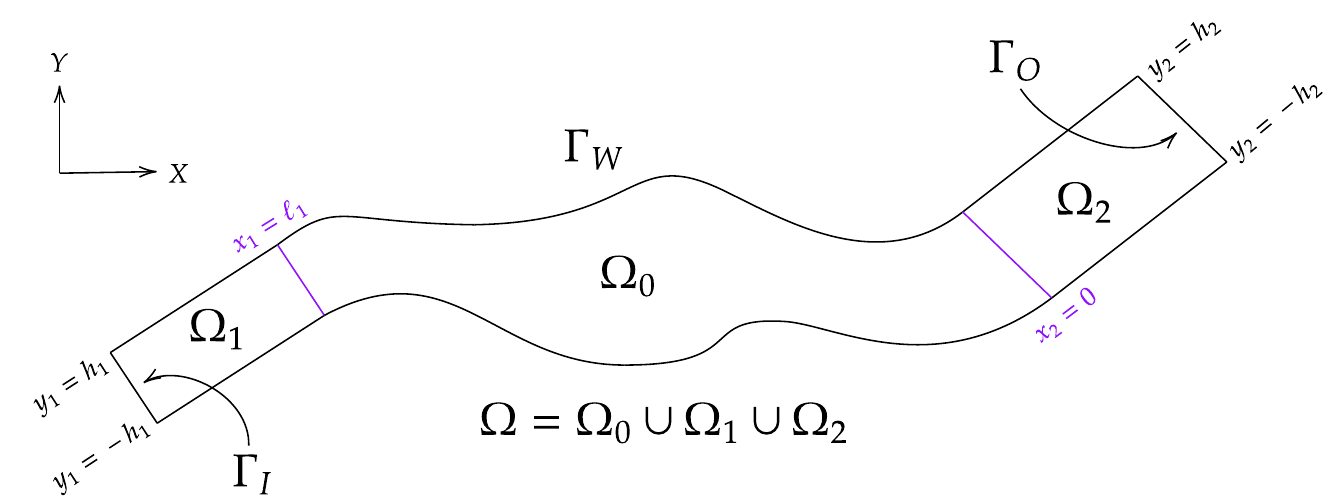}
	\end{center}
	\vspace*{-5mm}
	\caption{Representation of an admissible domain $\Omega \subset \mathbb{R}^{2}$.}\label{dom3}
\end{figure}

\begin{definition} \label{addomain3}
	An open bounded set $\Omega \subset \mathbb{R}^{3}$ will be called \textbf{admissible} if $\partial \Omega$ is piecewise of class $\mathcal{C}^{2}$, $\Omega$ is simply connected and $\Omega$ is the union of three disjoint subsets as follows:
	\begin{itemize}
		\item [(1)] In some coordinate system, $\Omega_{1} \doteq \Theta_{1} \times (0,\ell_{1})$, for some $\ell_{1}>0$ and some open, bounded and planar domain $\Theta_{1} \subset \mathbb{R}^{2}$ having boundary of class $\mathcal{C}^{2}$. Therefore, $\Omega_{1}$ is a cylinder of length $\ell_{1}$ and fixed cross-section, given by the planar set $\Theta_{1}$;
		\item [(2)] In some coordinate system, $\Omega_{2} \doteq \Theta_{2} \times (0,\ell_{2})$, for some $\ell_{2}>0$ and some open, bounded and planar domain $\Theta_{2} \subset \mathbb{R}^{2}$ having boundary of class $\mathcal{C}^{2}$. Therefore, $\Omega_{2}$ is a cylinder of length $\ell_{2}$ and fixed cross-section, given by the planar set $\Theta_{2}$;
		\item [(3)] $\Omega_{0} \doteq \Omega \setminus (\Omega_{1} \cup \Omega_{2})$  (note that $\Omega_{0}$ is not open).
	\end{itemize}
	The boundary of $\Omega$ is decomposed as $ \partial \Omega = \Gamma_{I} \cup \Gamma_{W} \cup \Gamma_{O}$, where
	\begin{equation}\label{boundaryomega3d}
	\begin{aligned}
	& \Gamma_{I} \doteq \Theta_{1} \times \{0\} \quad \text{(in the coordinate system defining $\Omega_{1}$)}  \,, \\[5pt]
	& \Gamma_{O} \doteq \Theta_{2} \times \{\ell_{2}\} \quad \text{(in the coordinate system defining $\Omega_{2}$)} \, ,
	\end{aligned}
\end{equation}	
	and $\Gamma_{W} \subset \mathbb{R}^{3}$ represents the $\mathcal{C}^{2}$-surface connecting $\Gamma_{I}$ with $\Gamma_{O}$.
\end{definition}
\begin{figure}[H]
	\begin{center}
		\includegraphics[scale=0.65]{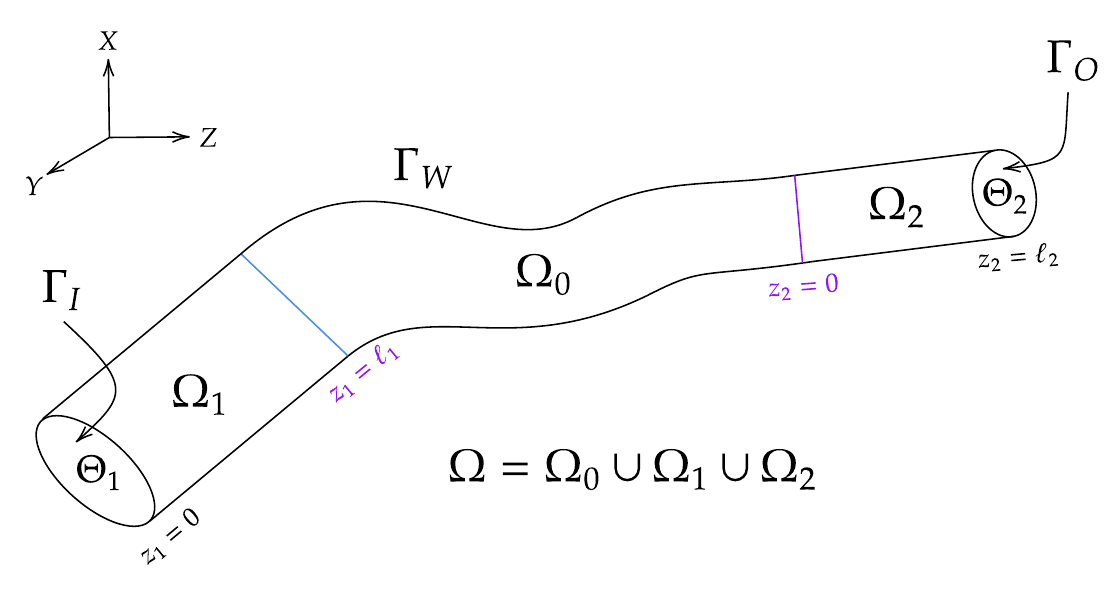}
	\end{center}
	\vspace*{-6mm}
	\caption{Representation of an admissible domain $\Omega \subset \mathbb{R}^{3}$.}\label{dom33d}
\end{figure}
\noindent
In few words, an admissible domain $\Omega \subset \mathbb{R}^{N}$, $N \in \{2,3\}$, in the sense of Definitions \ref{addomain2}-\ref{addomain3}, is a truncation (orthogonal to the symmetry axis of the inlet and outlet) of the domain considered in the celebrated \textbf{Leray problem}, as described in \cite[Definition 1.1]{amick1977steady} or \cite[Chapter III, Section 1.4]{galdi2008hemodynamical} (see also the chapter by Pileckas in \cite{pileckas2007navier}). Figure \ref{dom3} illustrates an example of such admissible planar channel, while Figure \ref{dom33d} depicts an admissible three-dimensional pipe.

Let $\Omega \subset \mathbb{R}^{N}$, $N \in \{2,3\}$, be an admissible domain in the sense of Definitions \ref{addomain2}-\ref{addomain3}. The outward unit normal to $\partial \Omega$ is denoted by $\bs{n} \in \mathbb{R}^{N}$. Henceforth we will refer to $\Gamma_{I}$ and $\Gamma_{O}$ in \eqref{boundaryomega1}-\eqref{boundaryomega3d} as the \textit{inlet} and \textit{outlet} of $\Omega$, respectively, while $\Gamma_{W}$ includes all the \textit{physical walls} of $\Omega$. 
For a fluid with viscosity $\nu > 0$, an inlet velocity field $\bs{g}_{*} : \Gamma_{I} \longrightarrow \mathbb{R}^{N}$ and a function $\sigma_{*} : \Gamma_{O} \longrightarrow \mathbb{R}$, by \textit{reference flow} we mean a pair $(\bs{W}_{\! \! *}, \Pi_{*}) \in W^{2,2}(\Omega;\mathbb{R}^N) \times W^{1,2}(\Omega;\mathbb{R})$ satisfying the following system in $\Omega$:
\begin{equation}\label{stokesintro}
	\left\{
	\begin{aligned}
		& \nabla\cdot \bs{W}_{\! \! *} = 0 \ \ \mbox{ in } \ \ \Omega \, , \\[5pt]
		& \bs{W}_{\! \! *}=\bs{g}_{*} \ \ \mbox{ on } \ \ \Gamma_{I}  \, , \quad \bs{W}_{\! \! *}=\bs{0} \ \ \mbox{ on } \ \ \Gamma_{W} \, , \\[5pt]
		& \nu \dfrac{\partial \bs{W}_{\! \! *}}{\partial \bs{n}} - \Pi_{*} \, \bs{n} = \sigma_{*} \, \bs{n} \ \ \mbox{ on } \ \ \Gamma_{O}  \, ,
	\end{aligned}
	\right.
\end{equation}
whose existence will be proved in Lemma \ref{lemma0}. In the present article we analyze the steady motion of a viscous incompressible fluid (identified by a viscosity equal to $\nu$) along $\Omega$, which is characterized by its velocity vector field $\bs{u} : \Omega \longrightarrow \mathbb{R}^N$ and its scalar pressure $p : \Omega \longrightarrow \mathbb{R}$, under the action of an external force $\bs{f} : \Omega \longrightarrow \mathbb{R}^N$. Such stationary motion will be modeled through the following boundary-value problem (with mixed boundary conditions) associated to the steady-state Navier-Stokes equations in $\Omega$ 
\begin{equation}\label{nsstokes0}
	\left\{
	\begin{aligned}
		& -\nu\Delta \bs{u}+(\bs{u}\cdot\nabla)\bs{u}+\nabla p=\bs{f} \, , \quad  \nabla\cdot \bs{u}=0 \ \ \mbox{ in } \ \ \Omega \, , \\[5pt]
		& \bs{u}=\bs{g}_{*} \ \ \mbox{ on } \ \ \Gamma_{I} \, , \\[5pt]
		& \bs{u}=\bs{0} \ \ \mbox{ on } \ \ \Gamma_{W} \, , \\[5pt]
		& \nu \dfrac{\partial \bs{u}}{\partial \bs{n}} - p \, \bs{n} + \dfrac{1}{2} [\bs{u} \cdot \bs{n}]^{-}(\bs{u} - \bs{W}_{\! \! *}) = \sigma_{*} \, \bs{n} \ \ \mbox{ on } \ \ \Gamma_{O} \, .
	\end{aligned}
	\right.
\end{equation}
While \eqref{nsstokes0}$_{2}$ prescribes the inlet velocity on $\Gamma_{I}$ and \eqref{nsstokes0}$_{3}$ describes the usual no-slip boundary condition on the physical walls $\Gamma_{W}$, identity \eqref{nsstokes0}$_{4}$ dictates that the fluid motion is subject to a generalized \textit{directional do-nothing} (DDN) boundary condition on $\Gamma_{O}$, see \cite{braack2014directional}, where 
\begin{equation} \label{plusminus}
[z]^{\pm} \doteq \dfrac{|z | \pm z}{2} \qquad \forall z \in \mathbb{R} \, .
\end{equation}

In order to properly motivate the choice of the boundary condition \eqref{nsstokes0}$_{4}$ on $\Gamma_{O}$, recall that it stems from a modification of the standard \textit{constant traction} or \textit{do-nothing} boundary condition, which reads
\begin{equation}\label{ctbc}
\nu \dfrac{\partial \bs{u}}{\partial \bs{n}} - p \, \bs{n} = \sigma_{*} \, \bs{n} \ \ \mbox{ on } \ \ \Gamma_{O} \, .
\end{equation}
From a mechanical point of view, the use of \eqref{ctbc} can be justified because the fluid flow, described on the inlet $\Gamma_{I}$ by the function $\bs{g}_{*}$, is altered along the channel or pipe, forbidding us to know in advance its behavior on the outlet. Hence, \eqref{ctbc} describes how the surroundings act on the fluid at the boundary $\Gamma_{O}$. For example, a pressure measured over $\Gamma_O$ can be imposed by a scalar function $\sigma_{*} \not= 0$, while for a traction-free outlet, we have $\sigma_{*} =0$, corresponding to the classical do-nothing condition. In turn, from a mathematical perspective, the boundary condition \eqref{ctbc} arises naturally while studying the weak formulation associated with the Navier-Stokes equations \eqref{nsstokes0}$_1$; indeed, suppose that $(\bs{u}, p) \in \mathcal{C}^{2}(\overline{\Omega}; \mathbb{R}^{N}) \times \mathcal{C}^{1}(\overline{\Omega}; \mathbb{R})$ is a classical solution to \eqref{nsstokes0}$_{1,2,3}$. Given any divergence-free vector field $\bs{\varphi} \in W^{1,2}(\Omega; \mathbb{R}^{N})$ that vanishes on $\Gamma_{I} \cup \Gamma_{W}$, upon integration by parts in $\Omega$ we obtain 
\begin{equation} \label{eqintro2}
\int_{\Omega} \left[ -\nu\Delta \bs{u}+(\bs{u}\cdot\nabla)\bs{u}+\nabla p  \right] \cdot \bs{\varphi} = \nu \int_{\Omega} \nabla \bs{u} : \nabla \bs{\varphi} + \int_{\Omega} (\bs{u}\cdot\nabla)\bs{u} \cdot \bs{\varphi} - \int_{\Gamma_{O}} \left( \nu \dfrac{\partial \bs{u}}{\partial \bs{n}} - p \, \bs{n} \right) \cdot \bs{\varphi} \, ,
\end{equation}
which puts \eqref{ctbc} as the foremost intuitive alternative to get rid of the boundary integral on the right-hand side of \eqref{eqintro2}. This is the reason why, since the formulation by Gresho \cite{gresho1991some} in 1991, the boundary condition \eqref{ctbc} has been widely employed in Computational Fluid Dynamics \cite{braack2014directional, heywood1996artificial, john2002higher, lanzendorfer2020multiple, rannacher2012short}. Nevertheless, it is well-known that \eqref{ctbc} does not allow for a control of the kinetic energy of the fluid flow, as the possibility of a backwards flow coming into $\Omega$ from $\Gamma_{O}$ is not excluded. Indeed, as asserted by Kra\v{c}mar \& Neustupa \cite{kravcmar2001weak, kravcmar2018modeling}, the transfer of kinetic energy from $\Gamma_{O}$ to $\Omega$ is due to the convective (non-linear) term $(\bs{u}\cdot\nabla)\bs{u}$, as the following integration by parts shows: 
\begin{equation} 
\int_\Omega (\bs{u}\cdot\nabla)\bs{u}\cdot \bs{u}=\dfrac{1}{2}\int_{\Gamma_I}|\bs{g}_*|^2(\bs{g}_*\cdot \bs{n})+\dfrac{1}{2}\int_{\Gamma_O}|\bs{u}|^2(\bs{u}\cdot \bs{n}) \, .
\label{introeq1}
\end{equation}
While the boundary integral on $\Gamma_{I}$ is completely controlled by the datum $\bs{g}_{*}$, there is no obvious way to dominate the integral over $\Gamma_O$ by other terms in the energy estimates. Recalling \eqref{plusminus}, we rewrite it as
\begin{equation} 
\dfrac{1}{2}\int_{\Gamma_O}|\bs{u}|^2(\bs{u}\cdot \bs{n})= \dfrac{1}{2}\int_{\Gamma_O}|\bs{u}|^2\left([\bs{u}\cdot \bs{n}]^+-[\bs{u}\cdot\bs{n}]^-\right) \, ,
\label{intbmn}
\end{equation}
which further reveals the possibility of backflow on $\Gamma_O$ whenever $[\bs{u}\cdot\bs{n}]^-\neq 0$ and suggests that the outflow boundary condition should account for this feature of the flow. Mathematically, as explained by Galdi in \cite[Chapter III]{galdi2008hemodynamical}, this translates into the fact that, up to nowadays, existence of solutions for the Navier-Stokes equations \eqref{nsstokes0}$_{1}$ with boundary conditions \eqref{nsstokes0}$_{2}$-\eqref{nsstokes0}$_{3}$-\eqref{ctbc} is known to hold only under a smallness assumption on the data of the problem (inlet velocity and external force). 

Based on \eqref{intbmn} and the preceding  considerations, one could justify an energy-stable boundary condition by modifying \eqref{ctbc} through the addition of the term $\frac 12 [\bs{u}\cdot\bs{n}]^-\bs{u}$. However, \eqref{nsstokes0}$_4$ also involves a reference velocity $\bs{W}_{\! \! *}$ whose choice has been rigorously  justified and further supported by numerical experiments in the works of Bruneau \& Fabrie \cite{bruneau2000boundary,BruneauFabrie1994,bruneau1996new}. Note that \eqref{introeq1} is formally obtained from \eqref{eqintro2}, replacing the test function $\bs{\varphi}$ with $\bs{u}$, but when $\bs{g}_{*}\not=0$ the boundary term $\int_{\Gamma_{I}} ( \nu \dfrac{\partial \bs{u}}{\partial \bs{n}} - p \, \bs{n} ) \cdot \bs{g}_{*} $ will appear in the energy equation, while the traction $\nu \dfrac{\partial \bs{u}}{\partial \bs{n}} - p \, \bs{n} $ on $\Gamma_I$  cannot be suitably bounded. 
The reference flow $(\bs{W}_{\! \! *}, \Pi_{*})$, which is not unique \cite{bruneau2000boundary}, incorporates the prescribed inlet velocity and outlet traction, providing a baseline representation of the fluid motion. Any deviation of the total flow from this reference is then represented  by the correction field $(\bs{v},q) := (\bs{u} - \bs{W}_{\! \! *}, p - \Pi_{*})$, whose velocity satisfies $\bs{v} = \bs{0}$ on $\Gamma_I$ and yields a contribution  for the kinetic energy balance in the following way:
\begin{equation}
\int_\Omega (\bs{u}\cdot\nabla)\bs{u}\cdot \bs{v} = \dfrac{1}{2}\int_{\Gamma_O}|\bs{v}|^2\left([\bs{u}\cdot \bs{n}]^+-[\bs{u}\cdot\bs{n}]^-\right) + \int_\Omega (\bs{u}\cdot\nabla) \bs{W}_{\! \! *} \cdot \bs{v}.
\label{ennovo}
\end{equation}
Physically, the term $(\bs{u}\cdot\nabla) \bs{W}_{\! \! *}$ represents the advection of the reference flow by the total velocity field $\bs{u}$ along the domain. The correction $\bs{v} = \bs{u} - \bs{W}_{\! \! *}$ thus allows the actual velocity field to deviate from the baseline field $\bs{W}_{\! \! *}$. Comparing \eqref{ennovo} with \eqref{introeq1}, the integral over $\Gamma_O$ now suggests that the outlet boundary condition should account for the convective transport of the velocity correction across the boundary by the total flow, through the term $\frac 12 [\bs{u}\cdot \bs{n}]^- \bs{v} $, which is activated only when backflow occurs.

Hence, a feasible alternative to the standard outflow condition \eqref{ctbc}, which overcomes the difficulties presented above, is given by the DDN-type boundary condition \eqref{nsstokes0}$_4$, as it has the advantage of being energy-stable and ensuring the uniqueness and stability of the rest state \cite{braack2014directional}. We stress out that both boundary conditions, \eqref{nsstokes0}$_4$ and \eqref{ctbc}, coincide in the case of pure outflow on $\Gamma_{O}$, namely, when $\bs{u} \cdot \bs{n} \geq 0$ on $\Gamma_{O}$. A particular example is given by the Hagen-Poiseuille flow, see Remarks \ref{remexact1}-\ref{remexact2} for further details. In Remark \ref{adim} we present the dimensionless form of \eqref{nsstokes0}$_4$ in terms of the Reynolds number. Computationally, this energy-stable boundary condition has been shown to contribute greatly in the prevention of spurious backflow at outflow boundary portions \cite{bruneau2000boundary,BCBBG2018,DONG2015300}. Mathematically, the fundamental goal of the present paper is to show that it allows for the unrestricted solvability (in a sense to be made precise later) of the system \eqref{nsstokes0}, in particular, without any assumption on the size of the  inlet velocity or the size of the corresponding flow rate, which is necessary in the case of pure Dirichlet boundary conditions \cite{galdi2011introduction,korobkov2015solution}. A significant advantage of the directional boundary condition \eqref{nsstokes0}$_4$ over other energy-stable outflow conditions in pipe systems, including boundary conditions formulated in terms of the Bernoulli pressure \cite{Baranovskii2020OptimalBC,FMN2007}, is that \eqref{nsstokes0}$_4$ is satisfied by Poiseuille flow. Note that the use of the Bernoulli pressure  in the Neumann boundary condition leads to a nonstandard formulation for the Navier-Stokes equations and the well-posedness results in \cite{FMN2007,Baranovskii2020OptimalBC} are established under homogeneous Dirichlet boundary conditions.

The method we employ to establish the existence of a weak solution of \eqref{nsstokes0} in Theorem \ref{epslevel}, the main result of the paper, is based on a contradiction argument that traces back to the seminal work of Leray \cite{leray1933etude}. The idea, which uses the Leray-Schauder Principle, has subsequently been extended and adapted in numerous mathematical studies concerning the steady-state Navier-Stokes equations, including the long-standing \textit{flux problem} solved in 2015 by Korobkov, Pileckas \& Russo \cite{korobkov2015solution} (see also their recent book \cite{korobkov2024steady}).  Here, we shall employ similar arguments to obtain uniform bounds on solutions by assuming the contrary, rescaling the equations, and then finding a solution to the steady Euler equations in $\Omega$ through a limiting process. A crucial ingredient is the fact that, as a consequence of the Bernoulli law, the pressure corresponding to the weak solution of the Euler equations takes a constant value on $\Gamma_I \cup \Gamma_W$, while it is a non-positive function on $\Gamma_O$. In the final step, the contradiction is derived by testing the Euler equations with $\bs{W}_{\! \! *}$, the extension constructed in Lemma \ref{lemma0}, followed by an integration by parts, where the sign of the flux of $\bs{g}_{*}$ across the inlet $\Gamma_I$, along with the specific structure of the energy-stable outflow boundary condition \eqref{nsstokes0}$_{4}$, play a key role.

The paper is organized as follows. Preliminary results are laid out in Section \ref{prelresults}, namely: the construction of the reference flow $(\bs{W}_{\! \! *}, \Pi_{*})$ in Lemma \ref{lemma0}, and a result concerning the inversion of the divergence operator on $L^{2}(\Omega;\mathbb{R})$ (Lemma \ref{bogtype}). Section \ref{secunres} is devoted to the regularity (Theorem \ref{press}) and existence (Theorem \ref{epslevel}) of weak solutions to system \eqref{nsstokes0}, without any restriction on the size of the data of the problem. Unique solvability of the boundary-value problem \eqref{nsstokes0} is established in Section \ref{solunica}, see Theorem \ref{theounique1}, under a suitable restriction on the magnitude of the data. 
\begin{remark}\label{adim}
	Another canonical form of the Navier-Stokes equations \eqref{nsstokes0} can be obtained by introducing the Reynolds number $\text{Re} = LU/\nu$, where $L$ and $U$ are characteristic length and velocity scales, respectively.  
	In this case the characteristic velocity is given by the bulk inflow velocity, namely,
	$$
	U = \frac{1}{|\Gamma_I|} \int_{\Gamma_I}|\bs{g}_* \cdot \bs{n}|,
	$$
	while the inlet diameter (or another relevant characteristic length of the fluid domain) determines $L$. Defining the dimensionless quantities 
	$$
	\bs{x}'=\frac{\bs{x}}{L}, \quad \bs{u}'=\frac{\bs{u}}{U}, \quad p'=\frac{p}{U^2}, \quad \bs{f}'=\frac{L}{U^2} \bs{f},  \quad \bs{g}_*'=\frac{\bs{g}_*}{U},  \quad \bs{W}'_{\! \! *} = \frac{\bs{W}_{\! \! *}}{U}, \quad  \sigma'_{*}  = \frac{1}{U^2} \sigma_{*},\quad	\bs{n}'(\bs{x}') = \bs{n}(\bs{x}) \, ,
	$$
	and introducing the rescaled domains
	$$
	\Omega' = \{ \bs{x}' = \bs{x} / L \mid \bs{x} \in \Omega\} \, ,  \qquad \Gamma'_{I}  = \{\bs{x}' = \bs{x} / L \mid \bs{x} \in \Gamma_{I}  \},  \quad (\text{similarly for} \ \ \Gamma_W' \ \text{ and } \ \Gamma'_{O}) \, ,
	$$
	we obtain \eqref{nsstokes0} in non-dimensional form: 
	\begin{equation}\label{NSad}
		\left\{
		\begin{aligned}
			&- \frac{1}{\text{Re}} \Delta_{\bs{x}'} \bs{u}' + (\bs{u}' \cdot \nabla_{\bs{x}'})\bs{u}' + \nabla_{\bs{x}'} p' = \bs{f}' \, , \qquad \nabla_{\bs{x}'} \cdot \bs{u}' = 0 \qquad \text{in} \ \ \ \Omega' \\[5pt]
			& \bs{u}'=\bs{g}'_{*} \ \ \mbox{ on } \ \ \Gamma'_{I} \, , \\[5pt]
			& \bs{u}'=\bs{0} \ \ \mbox{ on } \ \ \Gamma'_{W} \, , \\[5pt]
			& \frac{1}{\text{Re}} \dfrac{\partial \bs{u}'}{\partial \bs{n}'} - p' \, \bs{n}' + \dfrac{1}{2} [\bs{u}' \cdot \bs{n}']^{-}(\bs{u}' - \bs{W}'_{\! \! *}) = \sigma'_{*} \, \bs{n}' \ \ \mbox{ on } \ \ \Gamma'_{O}.
		\end{aligned}
		\right.
	\end{equation}
	Observe that this system retains the same structure as \eqref{nsstokes0}, with $\nu$ replaced by $1/\text{Re}$. Concerning the dimensionless form of the DDN boundary condition, we refer to \cite{bruneau2000boundary,bruneau1996new} for further details.
\end{remark}

\textbf{Notation:} Given any domain $D \subset \mathbb{R}^{N}$, $N \in \{2,3\}$, an integer $k \geq 1$ and an exponent $q \in [1, \infty]$, Lebesgue and Sobolev spaces will be denoted, respectively, by $L^q(D;{\mathbb R}^N)$ and $W^{k,q}(D;{\mathbb R}^N)$, and we write $\| \cdot \|_{L^q(D)}$ and $\| \cdot \|_{W^{k,q}(D)}$ for the corresponding norms. Given any boundary portion $\Gamma \subset \partial D$, Lebesgue and trace spaces will be denoted by $L^q(\Gamma;{\mathbb R}^N)$ and $W^{k - 1/q,q}(\Gamma;{\mathbb R}^N)$, respectively, with similar notations for their norms. Throughout the paper, otherwise specified, $C > 0$ will always denote a generic constant that depends exclusively on $\Omega$ and $\nu$, but that may change from line to line.

\section{Preliminary results} \label{prelresults}
Let $\Omega \subset \mathbb{R}^{N}$, $N \in \{2,3\}$, be an admissible domain occupied by an incompressible fluid with viscosity $\nu$. Since the portions $\Omega_{1}$ and $\Omega_{2}$ are given by rectangles ($N=2$) or cylinders ($N=3$), we can define a family of \textit{fully developed flows}, categorized by the fact that the only nonzero component of the associated velocity field is directed along the axis of $\Omega_{1}$ and $\Omega_{2}$. Given any $\Phi \geq 0$, consider the following:

\medskip
$\text{(I)}$ For $N=2$ and $j \in \{1,2\}$, introduce the horizontal rectangle $\mathcal{P}_{j} \doteq (0,\ell_{j}) \times (-h_{j},h_{j})$. The steady-state \textit{Hagen-Poiseuille flow} in the domain $\mathcal{P}_{j}$, with  flow rate $\Phi$, is characterized by the velocity field $\bs{U}^{(j)}_{\! \Phi} : \overline{\mathcal{P}_{j}} \longrightarrow \mathbb{R}^2$ and scalar pressure $Q^{(j)}_{\Phi} : \overline{\mathcal{P}_{j}} \longrightarrow \mathbb{R}$ defined as
\begin{equation}\label{poi1}
	\bs{U}^{(j)}_{\! \Phi}(x,y) \doteq \left(\dfrac{3\Phi}{4 h_{j}^{3}}(h_{j}^2 - y^2)\, , \, 0 \right)  \qquad \text{and} \qquad Q^{(j)}_{\Phi}(x,y) \doteq -\dfrac{3 \nu \Phi}{2 h_{j}^{3}} \left( x - \ell_{j} \right) \qquad \forall (x,y) \in \overline{\mathcal{P}_{j}} \, .
\end{equation}

$\text{(II)}$ For $N=3$ and $j \in \{1,2\}$, introduce the horizontal cylinder $\mathcal{P}_{j} \doteq \Theta_{j} \times (0,\ell_{j})$. Denote by $u_{j} \in \mathcal{C}^{2}(\overline{\Theta_{j}}; \mathbb{R})$ the unique classical solution to the following torsion problem in $\Theta_{j}$:
$$
-\Delta u_{j} = 1 \ \ \mbox{ in } \ \ \Theta_{j} \, , \qquad u_{j}=0 \ \ \mbox{ on } \ \ \partial \Theta_{j} \, ,
$$
and define
$$
\varrho_{j} \doteq \int_{\Theta_{j}} u_{j} = \int_{\Theta_{j}} | \nabla u_{j} |^{2} > 0 \, .
$$ 
By the Maximum Principle we have
\begin{equation} \label{maxprin}
	u_{j}(x,y) \geq 0 \qquad \forall (x,y) \in \overline{\Theta_{j}} \, .
\end{equation}
The steady-state Hagen-Poiseuille flow associated to $\mathcal{P}_{j}$, with flow rate $\Phi$, is characterized by the velocity field $\bs{U}^{(j)}_{\! \Phi} : \overline{\mathcal{P}_{j}} \longrightarrow \mathbb{R}^3$ and scalar pressure $Q^{(j)}_{\Phi} : \overline{\mathcal{P}_{j}} \longrightarrow \mathbb{R}$ defined as
\begin{equation}\label{poi23d}
	\bs{U}^{(j)}_{\! \Phi}(x,y,z) \doteq \dfrac{\Phi}{\varrho_{j}} (0,0, u_{j}(x,y) ) \qquad \text{and} \qquad Q^{(j)}_{\Phi}(x,y,z) \doteq -\dfrac{\nu \Phi}{\varrho_{j}} \left( z - \ell_{j} \right) \qquad \forall (x,y,z) \in \overline{\mathcal{P}_{j}} \, .
\end{equation}

Now, given the geometry of $\Omega_{j}$, there exists a point $\bs{x}_{j} \in \mathbb{R}^N$ and a rotation matrix $\mathcal{Q}_{j} \in \text{SO}(N)$ such that the set $\overline{\Omega_{j}}$ can be mapped onto $\overline{\mathcal{P}_{j}}$ through the following rigid transformation: 
$$
\mathbf{T}_{j}(\bs{x}) = \bs{x}_{j} + \mathcal{Q}_{j} \, \bs{x} \qquad \forall \bs{x} \in \overline{\Omega_{j}} \, .  
$$
Then, the steady-state Hagen-Poiseuille flow associated to $\Omega_{j}$, having flow rate $\Phi$, can be defined as
\begin{equation}\label{poi23dd}
	\bs{V}^{(j)}_{\! \Phi}(\bs{x}) \doteq \mathcal{Q}_{j}^{\top} \bs{U}^{(j)}_{\! \Phi}(\mathbf{T}_{j}(\bs{x})) \quad \text{and} \quad
	P_{\Phi}^{(j)}(\bs{x}) \doteq Q^{(j)}_{\! \Phi}(\mathbf{T}_{j}(\bs{x})) \qquad \forall \bs{x} \in \overline{\Omega_{j}} \, ,
\end{equation}
so that $(\bs{V}^{(j)}_{\! \Phi}, Q^{(j)}_{\! \Phi}) \in \mathcal{C}^{2}(\overline{\Omega_{2}}; \mathbb{R}^{N}) \times \mathcal{C}^{1}(\overline{\Omega_{2}}; \mathbb{R})$ solves in classical sense the following system:
\begin{equation}\label{flux03d}
	\left\{
	\begin{aligned}
		\begin{aligned}
			& -\nu\Delta \bs{V}^{(j)}_{\! \Phi} + \nabla P^{(j)}_{\Phi}  = \bs{0} \, , \quad  \nabla\cdot \bs{V}^{(j)}_{\! \Phi} = 0 \ \ \mbox{ in } \ \ \Omega_{j} \, , \\[5pt]
			& \bs{V}^{(j)}_{\! \Phi}=\bs{0} \ \ \mbox{ on } \ \ \Gamma_{W} \cap \partial \Omega_{j} \, , \\[5pt]
			& \int_{\Sigma} \bs{V}^{(j)}_{\! \Phi} \cdot \bs{n} = \Phi \quad \text{for any inner cross-section $\Sigma \subset \overline{\Omega_{j}}$ of $\Omega_{j}$} \, .
		\end{aligned}
	\end{aligned}
	\right.
\end{equation}

A smooth reference flow, capturing the inlet velocity on $\Gamma_{I}$ and subject to a do-nothing boundary condition on $\Gamma_{O}$, needs to be constructed. For this, recall the \textit{Lions-Magenes class} of inflows having a non-negative influx, defined as
$$
W^{3/2,2}_{+}(\Gamma_{I};\mathbb{R}^N) \doteq \left\{ \bs{g} \in W^{3/2,2}(\partial \Omega; \mathbb{R}^N) \ \Big| \ \bs{g} = \bs{0} \ \ \mbox{on} \ \ \partial\Omega \setminus \Gamma_{I} \, , \quad \int_{\Gamma_{I}} \bs{g} \cdot \bs{n} \leq 0 \, \right\} \, ,
$$
see \cite[Chapter VII]{dautray1999mathematical} for further details. In other words, $W^{3/2,2}_{+}(\Gamma_{I};\mathbb{R}^N)$ is the closed convex subset of $W^{3/2,2}(\partial \Omega; \mathbb{R}^N)$ comprising restrictions, to $\Gamma_{I}$, of vector fields having non-negative flux across $\Gamma_{I}$ and support contained in $\Gamma_{I}$. In fact, given $\bs{g}_{*} \in W^{3/2,2}_{+}(\Gamma_{I};\mathbb{R}^N)$, we define its \textit{flux across $\Gamma_{I}$} as
\begin{equation} \label{gflux}
	\Phi_{*} \doteq - \int_{\Gamma_{I}} \bs{g}_{*} \cdot \bs{n} \geq 0 \, .
\end{equation}

Additionally, given a bounded Lipschitz domain $D \subset \mathbb{R}^N$, consider the space of scalar functions with zero mean value in $D$, defined as
\begin{equation} \label{l02}
	L^2_0(D;\mathbb{R}) \doteq \left\{ g\in L^2(D;\mathbb{R}) \ \Big | \ \int_{D} g = 0 \right\} \, .
\end{equation}
It is well-known (see \cite[Section III.3]{galdi2011introduction} and Lemma \ref{bogtype} below) that this space plays a fundamental role in the inversion of the divergence operator on $L^2_0(D;\mathbb{R})$, and therefore, also for the obtainment of pressure estimates in the equations of hydrodynamics in general.

We then prove the following:

\begin{lemma}\label{lemma0}
Let $\Omega \subset \mathbb{R}^{N}$, $N \in \{2,3\}$, be an admissible domain. For any $\bs{g}_{*} \in W^{3/2,2}_{+}(\Gamma_{I};\mathbb{R}^N)$ and $\sigma_{*} \in W^{1/2,2}(\Gamma_{O};\mathbb{R})$, there exists a pair $(\bs{W}_{\! \! *}, \Pi_{*}) \in W^{2,2}(\Omega;\mathbb{R}^N) \times W^{1,2}(\Omega;\mathbb{R})$ such that
\begin{equation}\label{stopoi2}
	\left\{
\begin{aligned}
	& \nabla\cdot \bs{W}_{\! \! *} = 0 \ \ \mbox{ in } \ \ \Omega \, , \\[5pt]
	& \bs{W}_{\! \! *}=\bs{g}_{*} \ \ \mbox{ on } \ \ \Gamma_{I}  \, , \quad \bs{W}_{\! \! *}=\bs{0} \ \ \mbox{ on } \ \ \Gamma_{W} \, , \\[5pt]
	& \nu \dfrac{\partial \bs{W}_{\! \! *}}{\partial \bs{n}} - \Pi_{*} \, \bs{n} = \sigma_{*} \, \bs{n} \ \ \mbox{ on } \ \ \Gamma_{O} \, .
\end{aligned}
	\right.
\end{equation}
Moreover,  
\begin{equation}\label{flux}
\bs{W}_{\! \! *}=\bs{V}^{(2)}_{\! \Phi_{*}} \ \ \mbox{ on } \ \ \Gamma_{O} \, ,
\end{equation}
where $\bs{V}^{(2)}_{\! \Phi_{*}} \in \mathcal{C}^{2}(\overline{\Omega_{2}}; \mathbb{R}^{N})$ is the Hagen-Poiseuille flow defined in \eqref{poi23dd}, and the estimates
\begin{equation}\label{fluxestimate}
\| \bs{W}_{\! \! *} \|_{W^{2,2}(\Omega)} \leq C_{*} \| \bs{g}_{*} \|_{W^{3/2,2}(\Gamma_{I})} \qquad \text{and} \qquad \| \Pi_{*} \|_{W^{1,2}(\Omega)} \leq C_{*}  \| \sigma_{*} \|_{W^{1/2,2}(\Gamma_{O})}  \, ,
\end{equation}
hold for some constant $C_{*} > 0$ that depends exclusively on $\Omega$ and $\nu$.
\end{lemma}
\noindent
\begin{proof}	
In view of Definitions \ref{addomain2}-\ref{addomain3}, let us introduce (in local coordinates) the sets
$$
\begin{aligned}
& \Omega_{*} \doteq \Omega_{0} \cup \Omega_{1} \cup \left\lbrace (x_2,y_2) \in \Omega_{2} \ \Big\vert \  x_2 < \dfrac{\ell_{2}}{3} \, \right\rbrace \ \ \text{($N=2$)} \, ; \\[6pt]
& \Omega_{*} \doteq \Omega_{0} \cup \Omega_{1} \cup \left\lbrace (x_2,y_2,z_{2}) \in \Omega_{2} \ \Big\vert \  z_{2} < \dfrac{\ell_{2}}{3} \, \right\rbrace \ \ \text{($N=3$)} \, ,
\end{aligned}
$$
together with
$$
\Gamma_{*} \doteq \overline{\Omega} \cap \left\lbrace (x_2,y_2) \in \mathbb{R}^{2} \ \Big\vert \  x_{2} = \dfrac{\ell_{2}}{3} \, \right\rbrace \ \ \text{($N=2$)} \, ;\quad \Gamma_{*} \doteq \overline{\Omega} \cap \left\lbrace (x_2,y_2,z_{2}) \in \mathbb{R}^{3} \ \Big\vert \  z_{2} = \dfrac{\ell_{2}}{3} \, \right\rbrace \ \ \text{($N=3$)} \, .
$$	
Accordingly, let $(\bs{W}_{\! \! 0}, \Pi_{0}) \in W^{1,2}(\Omega_{*};\mathbb{R}^N) \times L_{0}^{2}(\Omega_{*};\mathbb{R})$ be the unique weak solution to the following Stokes system in $\Omega_{*}$:
\begin{equation}\label{stopoi2truncated}
	\left\{
	\begin{aligned}
		& -\nu\Delta \bs{W}_{\! 0} + \nabla \Pi_{0}  = \bs{0} \, , \quad  \nabla\cdot \bs{W}_{\! 0} = 0 \ \ \mbox{ in } \ \ \Omega_{*} \, , \\[5pt]
		& \bs{W}_{\! 0}=\bs{g}_{*} \ \ \mbox{ on } \ \ \Gamma_{I}  \, , \quad \bs{W}_{\! 0}=\bs{0} \ \ \mbox{ on } \ \ \Gamma_{W} \cap \partial \Omega_{*} \, , \\[5pt]
		& \bs{W}_{\! 0}=\bs{V}^{(2)}_{\! \Phi_{*}} \ \ \mbox{ on } \ \ \Gamma_{*} \, ,
	\end{aligned}
	\right.
\end{equation}
whose existence is guaranteed by standard methods (see, for example, \cite[Theorem IV.1.1]{galdi2011introduction}), owing to the compatibility condition emanating from \eqref{flux03d}$_3$ and \eqref{gflux}; in fact, we have
$$
\int_{\Gamma_{I}} \bs{g}_{*} \cdot \bs{n} + \int_{\Gamma_{*}} \bs{V}^{(2)}_{\! \Phi_{*}} \cdot \bs{n} = -\Phi_{*} + \Phi_{*} = 0 \, .
$$
Moreover, since the domains $\Omega_{1}$ and $\Omega_{2} \cap \Omega_{*}$ are convex polygons (when $N=2$) or cylinders (when $N=3$), and the lateral boundary $\Gamma_{W}$ is smooth, by merging the well-known regularity results for the steady-state Stokes equations under non-homogeneous Dirichlet boundary conditions (see \cite[Teorema, page 311]{cattabriga1961problema}) with a localization argument through a partition of unity (as in \cite[Theorem A.1]{concaenglish}) we can establish that $(\bs{W}_{\! \! 0}, \Pi_{0}) \in W^{2,2}(\Omega_{*};\mathbb{R}^N) \times W^{1,2}(\Omega_{*};\mathbb{R})$, together with the estimate
\begin{equation}\label{fluxestimate0}
	\| \bs{W}_{\! \! 0} \|_{W^{2,2}(\Omega_{*})} + \| \Pi_{0} \|_{W^{1,2}(\Omega_{*})} \leq C \left( \| \bs{g}_{*} \|_{W^{3/2,2}(\Gamma_{I})} + \| \bs{V}^{(2)}_{\! \Phi_{*}} \|_{W^{3/2,2}(\Gamma_{*})} \right) \, .
\end{equation} 

The next step consists in smoothly connecting the vector field $\bs{W}_{\! 0}$, defined in $\Omega_{*}$, with the velocity field $\bs{V}^{(2)}_{\! \Phi_{*}}$, defined in $\Omega_{2}$. For this, we further introduce the following subsets of $\Omega_{2}$:
$$
\begin{aligned}
& \Omega_{\sharp} \doteq \left\lbrace (x_2,y_2) \in \Omega_{2} \ \Big\vert \ 0 < x_{2} < \dfrac{\ell_{2}}{3} \, \right\rbrace \ \ \text{($N=2$)} \, ;\\[5pt]
& \Omega_{\sharp} \doteq \left\lbrace (x_2,y_2,z_{2}) \in \Omega_{2} \ \Big\vert \ 0 < z_{2} < \dfrac{\ell_{2}}{3} \, \right\rbrace \ \ \text{($N=3$)} \, ,
\end{aligned}
$$
together with
$$
\Gamma_{\sharp} \doteq \overline{\Omega} \cap \left\lbrace (x_2,y_2) \in \mathbb{R}^{2} \ \vert \  x_2 = 0 \, \right\rbrace \ \ \text{($N=2$)} \, ;\quad \Gamma_{\sharp} \doteq \overline{\Omega} \cap \left\lbrace (x_2,y_2,z_2) \in \mathbb{R}^{3} \ \vert \  z_2 = 0 \, \right\rbrace \ \ \text{($N=3$)} \, .
$$
Accordingly, take a cutoff function $\zeta \in \mathcal{C}^{\infty}_{0}(\mathbb{R}^N; [0,1])$ such that 
\begin{equation} \label{cutpro}
\zeta \equiv 1 \ \ \text{in a neighborhood of} \ \ \overline{\Omega_{0} \cup \Omega_{1}} \, ; \qquad \zeta \equiv 0 \ \ \text{in a neighborhood of} \ \ \overline{\Omega_{2} \setminus \overline{\Omega_{\sharp}}} \, ,
\end{equation}
see Figure \ref{cutoff} for an illustration of this configuration in the three-dimensional case.
\begin{figure}[H]
	\begin{center}
		\includegraphics[scale=0.8]{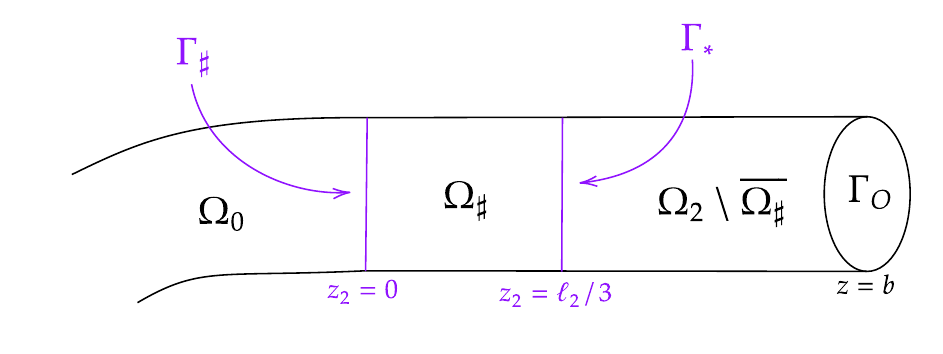}
	\end{center}
	\vspace*{-7mm}
	\caption{The domain $\Omega_{\sharp} \subset \Omega$ when $N=3$.}\label{cutoff}
\end{figure}
\noindent
We then introduce the vector field $\bs{Z}_{\! 0} \in W^{2,2}(\Omega_{\sharp}; \mathbb{R}^N)$ according to
$$
\bs{Z}_{\! 0} \doteq \zeta \bs{W}_{\! 0} + (1-\zeta) \bs{V}^{(2)}_{\! \Phi_{*}} \quad \text{in} \ \ \Omega_{\sharp} \, ,
$$
which satisfies
\begin{equation}\label{divz0}
\nabla \cdot \bs{Z}_{\! 0} = \nabla \zeta \cdot (\bs{W}_{\! 0} - \bs{V}^{(2)}_{\! \Phi_{*}}) \quad \text{in} \ \ \Omega_{\sharp} \, ,
\end{equation}
and therefore, $\nabla \cdot \bs{Z}_{\! 0} \in W_{0}^{1,2}(\Omega_{\sharp}; \mathbb{R})$, due to \eqref{cutpro}. Moreover, from \eqref{flux03d}$_3$-\eqref{gflux}-\eqref{stopoi2truncated}-\eqref{cutpro} and the Divergence Theorem, it clearly follows
$$
\int_{\Omega_{\sharp}} \nabla \cdot \bs{Z}_{\! 0} = \int_{\Gamma_{*}} \bs{Z}_{\! 0} \cdot \bs{n} + \int_{\Gamma_{\sharp}} \bs{Z}_{\! 0} \cdot \bs{n} = \int_{\Gamma_{*}} \bs{V}^{(2)}_{\! \Phi_{*}} \cdot \bs{n} + \int_{\Gamma_{\sharp}} \bs{W}_{\! 0} \cdot \bs{n} = \Phi_{*} - \Phi_{*} = 0 \, .
$$
Since $\Omega_{\sharp}$ is a Lipschitz domain, we can invoke \cite[Theorem III.3.3]{galdi2011introduction} to deduce the existence of another vector field $\bs{J}_{\! 0} \in W_{0}^{2,2}(\Omega_{\sharp}; \mathbb{R}^N)$ such that $\nabla \cdot \bs{J}_{\! 0} = - \nabla \cdot \bs{Z}_{\! 0}$ in $\Omega_{\sharp}$, and
\begin{equation}\label{fluxestimate1}
\| \bs{J}_{\! 0} \|_{W^{2,2}(\Omega_{\sharp})} \leq C \| \nabla \cdot \bs{Z}_{\! 0} \|_{W^{1,2}(\Omega_{\sharp})} \leq C \left( \| \bs{W}_{\! 0} \|_{W^{1,2}(\Omega_{*})} + \| \bs{V}^{(2)}_{\! \Phi_{*}} \|_{W^{1,2}(\Omega_{2})} \right) \, ,
\end{equation}
where the last inequality in \eqref{fluxestimate1} follows from \eqref{divz0}. We are then in position to define the vector field $\bs{W}_{\! \! *} \in W^{2,2}(\Omega; \mathbb{R}^N)$ by
\begin{equation}
	\bs{W}_{\! *} \doteq \begin{cases}
		\bs{W}_{\! 0} & \ \mathrm{ in } \ \ \Omega_{0} \cup \Omega_{1} \, ,  \\[4pt]
		\bs{J}_{\! 0} + \zeta \bs{W}_{\! 0} + (1-\zeta) \bs{V}^{(2)}_{\! \Phi_{*}} & \ \mathrm{ in }  \ \ \Omega_{\sharp} \, , \\[4pt]
		\bs{V}^{(2)}_{\! \Phi_{*}} & \ \mathrm{ in } \ \ \Omega_{2} \setminus \Omega_{\sharp} \, ,
	\end{cases}
\end{equation}
which, by the previous construction, satisfies \eqref{stopoi2}$_1$-\eqref{stopoi2}$_2$. Furthermore, since $\bs{W}_{\! *} = \bs{V}^{(2)}_{\! \Phi_{*}}$ in $\Omega_{2} \setminus \Omega_{\sharp}$, we then have the identity
\begin{equation} \label{normder0}
\dfrac{\partial \bs{W}_{\! \! *}}{\partial \bs{n}} = \bs{0} \ \ \mbox{ on } \ \ \Gamma_{O} \, .
\end{equation}
It then remains to associate a scalar pressure to the vector field $\bs{W}_{\! *}$, allowing, in combination with \eqref{normder0}, for the equality \eqref{stopoi2}$_3$ to be observed. Since $\sigma_{*} \in W^{1/2,2}(\Gamma_{O};\mathbb{R})$, there exists $\Pi_{*} \in W^{1,2}(\Omega;\mathbb{R})$ such that
$\Pi_{*} =  -\sigma_{*}$ on $\Gamma_{O}$, in the trace sense, together with the bound
\begin{equation}\label{fluxestimate2}
	\| \Pi_{*} \|_{W^{1,2}(\Omega)} \leq C \| \sigma_{*} \|_{W^{1/2,2}(\Gamma_{O})} \, .
\end{equation}
Now, noticing from \eqref{poi23dd}-\eqref{gflux} that
\begin{equation}\label{fluxestimate4}
\| \bs{V}^{(2)}_{\! \Phi_{*}} \|_{W^{2,2}(\Omega_{2})} + \| \bs{V}^{(2)}_{\! \Phi_{*}} \|_{W^{3/2,2}(\Gamma_{*})} \leq C \Phi_{*} \leq C \| \bs{g}_{*} \|_{W^{3/2,2}(\Gamma_{I})} \, ,
\end{equation}
the desired estimates \eqref{fluxestimate} are obtained as a consequence of \eqref{fluxestimate0}-\eqref{fluxestimate1}-\eqref{fluxestimate2}-\eqref{fluxestimate4}.
\end{proof}

Notice that the directional do-nothing boundary condition \eqref{nsstokes0}$_{4}$ forbids us from modifying the scalar pressure in \eqref{nsstokes0} by an additive constant. In particular, in contrast with the standard case of Dirichlet boundary conditions, we can no longer assume that the scalar pressure is an element of  $L_{0}^{2}(\Omega;\mathbb{R})$, see \eqref{l02}. Another essential preliminary result, therefore, concerns the inversion of the divergence operator on the space $L^{2}(\Omega;\mathbb{R})$, which is achieved by properly adapting the classical results involving $L_{0}^{2}(\Omega;\mathbb{R})$ (see, for example, \cite[Section III.3]{galdi2011introduction}) to the geometry of our setting. Inspired by \cite[Lemma 4.2]{korobkov2020solvability}, we now prove the following:
\begin{lemma} \label{bogtype}
Let $\Omega \subset \mathbb{R}^{N}$, $N \in \{2,3\}$, be an admissible domain.	Given $q \in L^{2}(\Omega;\mathbb{R})$, there exists a vector field $\bs{X} \in W^{1,2}(\Omega;\mathbb{R}^{N})$ such that
	\begin{equation} \label{vecje}
		\nabla \cdot \bs{X} = q \ \ \mbox{in} \ \ \Omega \, , \qquad \bs{X} = \bs{0} \ \ \mbox{on} \ \ \Gamma_{I} \cup \Gamma_{W} \qquad \text{and} \qquad \| \nabla \bs{X} \|_{L^{2}(\Omega)} \leq C_* \| q \|_{L^{2}(\Omega)} \, ,
	\end{equation}
	for some constant $C_{*} > 0$ that depends exclusively on $\Omega$.
\end{lemma}
\noindent
\begin{proof}
In what follows, $C > 0$ will always denote a generic constant that depends exclusively on $\Omega$, but that may change from line to line.

For $j \in \{1,2\}$, let $\bs{V}^{(j)}_{\! 1} \in W^{2,2}(\Omega_{j};\mathbb{R}^N)$ be as in \eqref{poi23dd}, with $\nu = \Phi =1$. Also, consider the unique weak solution $(\bs{Z}, \Pi) \in W^{1,2}(\Omega;\mathbb{R}^N) \times L_{0}^{2}(\Omega;\mathbb{R})$ to the following Stokes system in $\Omega$:
$$
	\left\{
	\begin{aligned}
		& -\Delta \bs{Z} + \nabla \Pi  = \bs{0} \, , \quad  \nabla\cdot \bs{Z} = 0 \ \ \mbox{ in } \ \ \Omega \, , \\[5pt]
		& \bs{Z}=\bs{V}^{(1)}_{\! 1} \ \ \mbox{ on } \ \ \Gamma_{I}  \, , \quad \bs{Z}=\bs{V}^{(2)}_{\! 1} \ \ \mbox{ on } \ \ \Gamma_{O} \, , \\[5pt]
		& \bs{Z}=\bs{0} \ \ \mbox{ on } \ \ \Gamma_{W} \, .
	\end{aligned}
	\right.
$$
Therefore, 	
\begin{equation} \label{hagen1}
	\int_{\Gamma_{O}} \bs{Z} \cdot \bs{n} = 1 \, .
\end{equation}
Now, let $B \subset \mathbb{R}^{N}$ be an open ball such that
$$
\overline{\Omega_{2}} \subsetneq B \qquad \text{and} \qquad \overline{\Omega_{1}} \subsetneq \mathbb{R}^N \setminus \overline{B} \, .
$$
Accordingly, take a cutoff function $\zeta \in \mathcal{C}^{\infty}_{0}(B; [0,1])$ such that $\zeta \equiv 1$ in $\overline{\Omega_{2}}$; in particular,
$$
\zeta = 0 \ \ \text{on} \ \ \Gamma_{I} \qquad \text{and} \qquad \zeta = 1 \ \ \text{on} \ \ \Gamma_{O} \, .
$$
Then, introduce the vector field
	\begin{equation} \label{bogtypevec}
		\bs{Q} \doteq \left( \int_{\Omega} q \right) \zeta \, \bs{Z} \qquad \text{in} \ \ \Omega \, ,
	\end{equation}
	which is clearly an element of $W^{1,2}(\Omega;\mathbb{R}^N)$ such that $\bs{Q} = \bs{0}$ on $\Gamma_{I} \cup \Gamma_{W}$, together with the estimate
	\begin{equation} \label{pointwise2}
		\| \nabla \bs{Q} \|_{L^{2}(\Omega)} \leq C \| q \|_{L^{2}(\Omega)} \, .
	\end{equation}
	On the other hand, from the Divergence Theorem and \eqref{hagen1} we obtain
	$$
	\begin{aligned}
		\int_{\Omega} \nabla \cdot \bs{Q} = \int_{\Gamma_{O}} \bs{Q} \cdot \bs{n} = \left( \int_{\Omega} q \right) \int_{\Gamma_{O}} \bs{Z} \cdot \bs{n} = \int_{\Omega} q \, ,
	\end{aligned}
	$$
	so that $q - \nabla \cdot \bs{Q} \in L^{2}_{0}(\Omega;\mathbb{R})$. Then, standard results such as \cite[Lemma 1]{pileckas1983spaces} (see also \cite[Lemma 1.8.13]{korobkov2024steady}) ensure the existence of another vector field $\bs{Y} \in W^{1,2}_{0}(\Omega;\mathbb{R}^N)$ such that
	\begin{equation} \label{bogtypevec2}
		\nabla \cdot \bs{Y} = q - \nabla \cdot \bs{Q} \ \ \text{in} \ \ \Omega \qquad \text{and} \qquad  \| \nabla \bs{Y} \|_{L^{2}(\Omega)} \leq C \| q - \nabla \cdot \bs{Q} \|_{L^{2}(\Omega)}  \, .
	\end{equation}
	We set $\bs{X} \doteq \bs{Y} + \bs{Q}$ which, in view of \eqref{pointwise2}-\eqref{bogtypevec2}, is an element of $W^{1,2}(\Omega;\mathbb{R}^{N})$ satisfying \eqref{vecje}.
\end{proof}

\section{Unrestricted solvability of the boundary-value problem} \label{secunres}
We introduce the function spaces (of vector fields) that will be employed hereafter:
\begin{equation}\label{funcspaces}
	\begin{aligned}
		& \mathcal{V}(\Omega) \doteq \left\lbrace \bs{v} \in W^{1,2}(\Omega;\mathbb{R}^N) \ | \ \nabla \cdot \bs{v}=0 \ \ \mbox{in} \ \ \Omega \, , \quad \bs{v} = \bs{0} \ \ \mbox{on} \ \ \Gamma_{W} \right\rbrace \, , \\[3pt]
		& \mathcal{H}_{*}(\Omega) \doteq \left\lbrace \bs{v} \in W^{1,2}(\Omega;\mathbb{R}^N) \ | \ \bs{v} = \bs{0} \ \ \mbox{on} \ \ \Gamma_{I} \cup \Gamma_{W} \right\rbrace , \\[3pt]
		& \mathcal{V}_{*}(\Omega) \doteq \left\lbrace \bs{v} \in W^{1,2}(\Omega;\mathbb{R}^N) \ | \ \nabla \cdot \bs{v}=0 \ \ \mbox{in} \ \ \Omega \, , \quad \bs{v} = \bs{0} \ \ \mbox{on} \ \ \Gamma_{I} \cup \Gamma_{W} \right\rbrace \, ,
	\end{aligned}
\end{equation}
which are Hilbert spaces if endowed with the Dirichlet scalar product of the gradients, denoted by
\begin{equation} \label{dirscal}
	[\bs{v}, \bs{w}]_{\Omega} \doteq \int_{\Omega} \nabla \bs{v} : \nabla \bs{w} \qquad \forall \bs{v},\bs{w} \in W^{1,2}(\Omega;\mathbb{R}^N) \, .
\end{equation}
We recall the definition of weak solutions for the velocity component of problem \eqref{nsstokes0}:
\begin{definition}\label{weaksolution}
Let $\Omega \subset \mathbb{R}^{N}$, $N \in \{2,3\}$, be an admissible domain. Let $\bs{f} \in L^2(\Omega;\mathbb{R}^N)$ be an external force, $\bs{g}_{*} \in W^{3/2,2}_{+}(\Gamma_{I};\mathbb{R}^N)$ an inflow, $\sigma_{*} \in W^{1/2,2}(\Gamma_{O};\mathbb{R})$ and $(\bs{W}_{\! \! *}, \Pi_{*}) \in W^{2,2}(\Omega;\mathbb{R}^N) \times W^{1,2}(\Omega;\mathbb{R})$ the reference flow from Lemma \ref{lemma0}. A vector field $\bs{u} \in \mathcal{V}(\Omega)$ is called a \textbf{weak solution} of problem \eqref{nsstokes0} if $\bs{u} - \bs{W}_{\! \! *} \in \mathcal{V}_{*}(\Omega)$ and
\begin{equation} \label{nstokesdebil1}
	\nu \int_{\Omega} \nabla \bs{u} : \nabla \bs{\varphi} + \int_{\Omega} (\bs{u} \cdot \nabla)\bs{u} \cdot \bs{\varphi} + \dfrac{1}{2} \int_{\Gamma_{O}} [\bs{u} \cdot \bs{n}]^{-}(\bs{u} - \bs{W}_{\! \! *}) \cdot \bs{\varphi}  = \int_{\Omega} \bs{f} \cdot \bs{\varphi} + \int_{\Gamma_{O}} \sigma_{*} (\bs{\varphi} \cdot \bs{n})  \, ,
\end{equation}
for every $\bs{\varphi} \in \mathcal{V}_{*}(\Omega)$.
\end{definition}

\begin{remark} \label{rem1}
We emphasize that the integrals over $\Gamma_{O}$ appearing in the weak formulation \eqref{nstokesdebil1} are bounded in the function spaces considered, as the trace embedding $W^{1,2}(\Omega;\mathbb{R}^N) \subset L^{4}(\partial \Omega;\mathbb{R}^N)$ (which holds for $N \in \{2,3\}$) yields
$$
\begin{aligned}
\left| \int_{\Gamma_{O}} [\bs{u} \cdot \bs{n}]^{-}(\bs{u} - \bs{W}_{\! \! *}) \cdot \bs{\varphi} \right| & \leq \| \bs{u} \|_{L^{4}(\Gamma_{O})} \| \bs{u} - \bs{W}_{\! \! *} \|_{L^{4}(\Gamma_{O})} \| \bs{\varphi} \|_{L^{2}(\Gamma_{O})} \\[3pt]
& \leq C_{*} \| \nabla \bs{u} \|_{L^{2}(\Omega)} \| \nabla (\bs{u} - \bs{W}_{\! \! *}) \|_{L^{2}(\Omega)} \| \nabla \bs{\varphi} \|_{L^2(\Omega)} \, ,
\end{aligned}
$$
for every $(\bs{u},\bs{\varphi}) \in \mathcal{V}(\Omega) \times \mathcal{V}_{*}(\Omega)$, for some constant $C_{*} >0$ that depends exclusively on $\Omega$. 
\end{remark}

In the next result we explain how a uniquely defined scalar pressure can be associated to the velocity component of any weak solution to problem \eqref{nsstokes0}. Moreover, we prove that the resulting pair is a strong solution of \eqref{nsstokes0}.

\begin{theorem} \label{press}
	Let $\Omega \subset \mathbb{R}^{N}$, $N \in \{2,3\}$, be an admissible domain. Given a force $\bs{f} \in L^2(\Omega;\mathbb{R}^N)$, $\bs{g}_{*} \in W^{3/2,2}_{+}(\Gamma_{I};\mathbb{R}^N)$ and $\sigma_{*} \in W^{1/2,2}(\Gamma_{O};\mathbb{R})$, suppose that $\bs{u}\in \mathcal{V}(\Omega)$ is a weak solution of problem \eqref{nsstokes0}, in the sense of Definition \ref{weaksolution}. Then, there exists a unique pressure $p \in L^{2}(\Omega;\mathbb{R})$ such that 
\begin{equation} \label{nstokesdebilpressure1}
\begin{aligned}
&	\nu \int_{\Omega} \nabla \bs{u} : \nabla \bs{\varphi} + \int_{\Omega} (\bs{u} \cdot \nabla)\bs{u} \cdot \bs{\varphi} - \int_{\Omega} p (\nabla \cdot \bs{\varphi}) + \dfrac{1}{2} \int_{\Gamma_{O}} [\bs{u} \cdot \bs{n}]^{-}(\bs{u} - \bs{W}_{\! \! *}) \cdot \bs{\varphi} \\[6pt]
& \hspace{-4mm} = \int_{\Omega} \bs{f} \cdot \bs{\varphi} + \int_{\Gamma_{O}} \sigma_{*} (\bs{\varphi} \cdot \bs{n})  \qquad \forall \bs{\varphi} \in \mathcal{H}_{*}(\Omega) \, .
\end{aligned}
\end{equation}	
Moreover, $(\bs{u},p) \in W^{2,2}(\Omega;\mathbb{R}^N) \times W^{1,2}(\Omega;\mathbb{R})$, so that the pair $(\bs{u},p)$ satisfies \eqref{nsstokes0}$_1$ in strong sense in $\Omega$. The boundary conditions for $\bs{u}$ in \eqref{nsstokes0}$_2$-\eqref{nsstokes0}$_3$ are verified in the sense of $W^{3/2,2}(\partial \Omega;\mathbb{R}^N)$, while the boundary condition \eqref{nsstokes0}$_4$ is verified in the sense of $W^{1/2,2}(\Gamma_{O};\mathbb{R}^N)$. Additionally, there holds the bound
\begin{equation} \label{stimacatta00}
	\begin{aligned}
		\| \bs{u} \|_{W^{2,3/2}(\Omega)} + \| p \|_{W^{1,3/2}(\Omega)} & \leq C_{*} \Big( \| \bs{f} \|_{L^{2}(\Omega)} + \| \nabla \bs{u} \|^{2}_{L^{2}(\Omega)} + \| \nabla \bs{W}_{\! \! *} \|^{2}_{L^{2}(\Omega)}  \\[6pt] 
		& \hspace{1.2cm} + \| \bs{g}_{*} \|_{W^{3/2,2}(\Gamma_{I})} + \| \sigma_{*} \|_{W^{1/2,2}(\Gamma_{O})} \Big) \, ,
	\end{aligned}
\end{equation}
for some constant $C_{*} > 0$ that depends exclusively on $\Omega$ and $\nu$.
\end{theorem}
\noindent
\begin{proof}
We follow closely the proof of \cite[Theorem III.5.3]{galdi2011introduction} (see also \cite[Lemma 2.1]{girault2012finite}). Denote by $\mathcal{H}^{-1}(\Omega)$ the dual space of $\mathcal{H}_{*}(\Omega)$ (see \eqref{funcspaces}$_2$) and by $\langle\cdot,\cdot\rangle_{\mathcal{H}_{*}(\Omega)}$ the
duality product
between $\mathcal{H}^{-1}(\Omega)$ and $\mathcal{H}_{*}(\Omega)$. The adjoint of the (strong) divergence operator $\text{div}_{\Omega} : \mathcal{H}_{*}(\Omega) \longrightarrow L^{2}(\Omega;\mathbb{R})$ is defined by the map $\text{div}^{*}_{\Omega} : L^{2}(\Omega;\mathbb{R}) \longrightarrow \mathcal{H}^{-1}(\Omega)$ such that
$$
\langle \text{div}^{*}_{\Omega}(q) , \bs{\varphi} \rangle_{\mathcal{H}_{*}(\Omega)} \doteq \int_{\Omega} q (\nabla \cdot \bs{\varphi}) \qquad \forall q \in L^{2}(\Omega;\mathbb{R}) \, , \quad \forall \bs{\varphi} \in \mathcal{H}_{*}(\Omega) \, .
$$
Therefore, the Closed Range Theorem of Banach \cite[Chapter VII]{yosida2012functional} can be applied to deduce that
$$
\text{Range}(\text{div}^{*}_{\Omega}) = \text{Ker}(\text{div}_{\Omega})^{\perp} = \mathcal{V}_{*}(\Omega)^{\perp} \doteq \{ \mathcal{F} \in
\mathcal{H}^{-1}(\Omega) \ | \ \langle \mathcal{F} , \varphi \rangle_{\mathcal{H}_{*}(\Omega)} = 0 \quad \forall \varphi \in \mathcal{V}_{*}(\Omega) \, \} \, .
$$
Given a weak solution $\bs{u}\in \mathcal{V}(\Omega)$ of \eqref{nsstokes0}, we define
$\mathcal{F}_{\bs{u}} \in \mathcal{H}^{-1}(\Omega)$ by
$$
\begin{aligned}
\mathcal{F}_{\bs{u}}(\varphi) &	\doteq \nu \int_{\Omega} \nabla \bs{u} : \nabla \bs{\varphi} + \int_{\Omega} (\bs{u} \cdot \nabla)\bs{u} \cdot \bs{\varphi} + \dfrac{1}{2} \int_{\Gamma_{O}} [\bs{u} \cdot \bs{n}]^{-}(\bs{u} - \bs{W}_{\! \! *}) \cdot \bs{\varphi} \\[6pt]
	& \hspace{4mm} - \int_{\Omega} \bs{f} \cdot \bs{\varphi} - \int_{\Gamma_{O}} \sigma_{*} (\bs{\varphi} \cdot \bs{n})  \qquad \forall \bs{\varphi} \in \mathcal{H}_{*}(\Omega) \, .
\end{aligned}
$$
Due to \eqref{nstokesdebil1}, $\mathcal{F}_{\bs{u}}\in
\text{Range}(\text{div}^{*}_{\Omega})$. This ensures the existence of a pressure $p \in L^{2}(\Omega;\mathbb{R})$ verifying \eqref{nstokesdebilpressure1}. Suppose there exists another function $q \in L^{2}(\Omega;\mathbb{R})$ satisfying \eqref{nstokesdebilpressure1}, thereby implying
\begin{equation} \label{bog001}
	\int_{\Omega} (p-q) (\nabla \cdot \bs{\varphi}) = 0 \qquad \forall \bs{\varphi} \in \mathcal{H}_{*}(\Omega) \, .
\end{equation}
Lemma \ref{bogtype} provides the existence of a vector field $\bs{X} \in \mathcal{H}_{*}(\Omega)$ such that $\nabla \cdot \bs{X}  = p-q$ in $\Omega$. Taking $\bs{\varphi} = \bs{X} $ in \eqref{bog001} ensures that $p=q$ almost everywhere in $\Omega$. 

Define $\bs{v} \doteq \bs{u} - \bs{W}_{\! \! *}$ and $q \doteq p - \Pi_{*}$, so that $(\bs{v},q) \in \mathcal{V}_{*}(\Omega) \times L^{2}(\Omega;\mathbb{R})$ satisfies the system
\begin{equation}\label{nsstokes0v}
	\left\{
	\begin{aligned}
		& -\nu\Delta \bs{v}+\nabla q=\bs{f} + \nu \Delta \bs{W}_{\! \! *} - \nabla \Pi_{*} - (\bs{v} \cdot \nabla)\bs{W}_{\! \! *} - (\bs{W}_{\! \! *} \cdot \nabla)\bs{v} - (\bs{W}_{\! \! *} \cdot \nabla)\bs{W}_{\! \! *} - (\bs{v}\cdot\nabla)\bs{v} \ \ \mbox{ in } \ \ \Omega \, , \\[5pt]
		& \nabla\cdot \bs{v}=0 \ \ \mbox{ in } \ \ \Omega \, , \\[5pt]
		& \bs{v}=\bs{0} \ \ \mbox{ on } \ \ \Gamma_{I} \cup  \Gamma_{W} \, , \\[5pt]
		& \nu \dfrac{\partial \bs{v}}{\partial \bs{n}} - q \, \bs{n} + \dfrac{1}{2} [(\bs{v} + \bs{W}_{\! \! *}) \cdot \bs{n}]^{-} \bs{v} = \bs{0} \ \ \mbox{ on } \ \ \Gamma_{O} 
	\end{aligned}
	\right.
\end{equation}
in weak form. Now, the Sobolev embedding $W^{1,2}(\Omega; \mathbb{R}^N) \subset L^{6}(\Omega; \mathbb{R}^N)$ implies $(\bs{v} \cdot \nabla)\bs{v} \in L^{3/2}(\Omega; \mathbb{R}^N)$, and as a consequence, that
$$
\bs{F} \doteq \bs{f} + \nu \Delta \bs{W}_{\! \! *} - \nabla \Pi_{*} - (\bs{v} \cdot \nabla)\bs{W}_{\! \! *} - (\bs{W}_{\! \! *} \cdot \nabla)\bs{v} - (\bs{W}_{\! \! *} \cdot \nabla)\bs{W}_{\! \! *} - (\bs{v} \cdot \nabla)\bs{v} \in L^{3/2}(\Omega; \mathbb{R}^N) \, .
$$ 
From the H\"older, Sobolev, Young and Poincaré inequalities we infer the estimate
\begin{equation}\label{stima0}
\begin{aligned}
	\| \bs{F} \|_{L^{3/2}(\Omega)} \leq C \left( \| \bs{f} \|_{L^{2}(\Omega)} + \| \bs{W}_{\! \! *} \|_{W^{2,2}(\Omega)} + \| \Pi_{*} \|_{W^{1,2}(\Omega)} + \| \nabla \bs{v} \|^{2}_{L^{2}(\Omega)} + \| \nabla \bs{W}_{\! \! *} \|^{2}_{L^{2}(\Omega)} \right) \, .
\end{aligned}
\end{equation}
Now, the restriction of $\bs{n}$ to $\Gamma_{O}$ is a constant vector, denoted by $\bs{n}_{*} \in \mathbb{R}^N$. We certainly have 
$$
[(\bs{v} + \bs{W}_{\! \! *}) \cdot \bs{n}_{*}]^{-} \in W^{1,2}(\Omega; \mathbb{R}) \, ,
$$
together with
\begin{equation}\label{stima1}
\| \nabla [(\bs{v} + \bs{W}_{\! \! *}) \cdot \bs{n}_{*}]^{-} \|_{L^{2}(\Omega)} \leq  \| \nabla \bs{v} \|_{L^{2}(\Omega)} + \| \nabla \bs{W}_{\! \! *} \|_{L^{2}(\Omega)} \, .
\end{equation}
By Sobolev embedding it also holds $[(\bs{v} + \bs{W}_{\! \! *}) \cdot \bs{n}_{*}]^{-} \bs{v} \in W^{1,3/2}(\Omega; \mathbb{R}^N)$. From \eqref{stima1} and the inequalities of H\"older, Sobolev and Young we obtain the bound
\begin{equation}\label{stima2000}
\begin{aligned}
\| \nabla ([(\bs{v} + \bs{W}_{\! \! *}) \cdot \bs{n}_{*}]^{-} \bs{v}) \|_{L^{3/2}(\Omega)} \leq C \left( \| \nabla \bs{v} \|^{2}_{L^{2}(\Omega)} + \| \nabla \bs{W}_{\! \! *} \|^{2}_{L^{2}(\Omega)} \right) \, .
\end{aligned}
\end{equation}
In particular, the previous observations imply that
$$
\bs{h} \doteq - \dfrac{1}{2} [(\bs{v} + \bs{W}_{\! \! *}) \cdot \bs{n}]^{-} \bs{v}  \in W^{\frac{1}{3},\frac{3}{2}}(\Gamma_{O}; \mathbb{R}^N) \, ,
$$
with
\begin{equation}\label{stima2}
\| \bs{h} \|_{W^{\frac{1}{3},\frac{3}{2}}(\Gamma_{O})} \leq C \| \nabla ([(\bs{v} + \bs{W}_{\! \! *}) \cdot \bs{n}]^{-} \bs{v}) \|_{L^{3/2}(\Omega)}  \leq C \left( \| \nabla \bs{v} \|^{2}_{L^{2}(\Omega)} + \| \nabla \bs{W}_{\! \! *} \|^{2}_{L^{2}(\Omega)} \right) \, ,
\end{equation}
as a consequence of the trace inequality and \eqref{stima2000}. In conclusion, the pair $(\bs{v},q) \in \mathcal{V}_{*}(\Omega) \times L^{2}(\Omega;\mathbb{R})$ can be considered to be a weak solution to the following Stokes system in $\Omega$:
\begin{equation}\label{stokes0}
	\left\{
	\begin{aligned}
		& -\nu\Delta \bs{v}+\nabla q=\bs{F} \, , \quad  \nabla\cdot \bs{v}=0 \ \ \mbox{ in } \ \ \Omega \, , \\[5pt]
		& \bs{v}=\bs{0} \ \ \mbox{ on } \ \ \Gamma_{I} \cup  \Gamma_{W} \, , \\[5pt]
		& \nu \dfrac{\partial \bs{v}}{\partial \bs{n}} - q \, \bs{n} = \bs{h} \ \ \mbox{ on } \ \ \Gamma_{O} \, .
	\end{aligned}
	\right.
	\end{equation}
Now, let $\bs{\chi} \in W^{1,2}(\Omega;\mathbb{R}^N)$ be the unique weak solution to the following Zaremba problem in $\Omega$:
\begin{equation}\label{mixedbc}
	\left\{
	\begin{aligned}
		& \Delta \bs{\chi}=\bs{0} \ \ \mbox{ in } \ \ \Omega \, , \\[5pt]
		& \bs{\chi}=\bs{0} \ \ \mbox{ on } \ \ \Gamma_{I} \cup  \Gamma_{W} \, , \\[5pt]
		& \dfrac{\partial \bs{\chi}}{\partial \bs{n}} = \dfrac{1}{\nu} \bs{h} \ \ \mbox{ on } \ \ \Gamma_{O} \, ,
	\end{aligned}
	\right.
\end{equation}
whose existence follows from the Lax-Milgram Theorem. Since the sets $\Omega_{1}$ and $\Omega_{2}$ correspond to rectangles ($N=2$) or cylinders ($N=3$), and the boundary portions $\Gamma_{W}$ and $\Gamma_{O}$ meet at a right angle, regularity results such as \cite[Chapter 4]{grisvard2011elliptic} (for $N=2$) and \cite[Corollary 8.3.2]{mazya2010elliptic} (for $N=3$; see also \cite{ott2013mixed} for a relevant discussion in both cases $N \in \{2,3\}$) ensure that $\bs{\chi} \in W^{2,3/2}(\Omega; \mathbb{R}^N)$, together with 
\begin{equation} \label{stimazaremba}
\| \bs{\chi} \|_{W^{2,3/2}(\Omega)} \leq C \| \bs{h} \|_{W^{\frac{1}{3},\frac{3}{2}}(\Gamma_{O})}  \, .
\end{equation}
We then set $\bs{z} \doteq \bs{v} - \bs{\chi}$ in $\Omega$, so that, due to \eqref{stokes0}-\eqref{mixedbc}, the pair $(\bs{z},q) \in \mathcal{V}_{*}(\Omega) \times L^{2}(\Omega;\mathbb{R})$ constitutes a weak solution to the following Stokes system in $\Omega$:
\begin{equation}\label{stokeszaremba}
	\left\{
	\begin{aligned}
		& -\nu\Delta \bs{z}+\nabla q=\bs{F} \, , \quad  \nabla\cdot \bs{z} = - \nabla \cdot \bs{\chi} \ \ \mbox{ in } \ \ \Omega \, , \\[5pt]
		& \bs{z}=\bs{0} \ \ \mbox{ on } \ \ \Gamma_{I} \cup  \Gamma_{W} \, , \\[5pt]
		& \nu \dfrac{\partial \bs{z}}{\partial \bs{n}} - q \, \bs{n} = \bs{0} \ \ \mbox{ on } \ \ \Gamma_{O} \, .
	\end{aligned}
	\right.
\end{equation}
We can then invoke \cite[Theorem 3.1]{neustupa2022maximum} and \cite[Corollary A.3]{benevs2014note} (when $N=2$ and $N=3$, respectively, see also \cite[Theorem A.1]{benevs2016solutions} for a closely related result when $N=2$; notice that the non-zero divergence condition in \eqref{stokeszaremba}$_1$ is irrelevant from a regularity perspective, as shown by \cite[Equation (A.10)]{benevs2016solutions} and \cite[Equation (A.6)]{benevs2014note}) to deduce that $(\bs{z},q) \in W^{2,3/2}(\Omega;\mathbb{R}^N) \times W^{1,3/2}(\Omega;\mathbb{R})$, together with
\begin{equation} \label{stima00}
\begin{aligned}
& \| \bs{z} \|_{W^{2,3/2}(\Omega)} + \| q \|_{W^{1,3/2}(\Omega)} \leq C \left( \| \bs{F} \|_{L^{3/2}(\Omega)} + \| \nabla \cdot \bs{\chi} \|_{W^{1,3/2}(\Omega)} \right) \\[6pt]
& \hspace{-4mm} \leq C \left( \| \bs{f} \|_{L^{2}(\Omega)} + \| \bs{W}_{\! \! *} \|_{W^{2,2}(\Omega)} + \| \Pi_{*} \|_{W^{1,2}(\Omega)} + \| \nabla \bs{v} \|^{2}_{L^{2}(\Omega)} + \| \nabla \bs{W}_{\! \! *} \|^{2}_{L^{2}(\Omega)} \right) \, ,
\end{aligned}
\end{equation}
	where the second inequality is obtained as a consequence of \eqref{stima0}-\eqref{stima2}-\eqref{stimazaremba}. Therefore,
\begin{equation} \label{stima000}
	\begin{aligned}
	\| \bs{v} \|_{W^{2,3/2}(\Omega)} + \| q \|_{W^{1,3/2}(\Omega)} & \leq C \Big( \| \bs{f} \|_{L^{2}(\Omega)} + \| \nabla \bs{v} \|^{2}_{L^{2}(\Omega)} + \| \nabla \bs{W}_{\! \! *} \|^{2}_{L^{2}(\Omega)}  \\[6pt] 
	& \hspace{1.2cm} + \| \bs{W}_{\! \! *} \|_{W^{2,2}(\Omega)} + \| \Pi_{*} \|_{W^{1,2}(\Omega)} \Big) \, .
\end{aligned}
\end{equation}
Recalling that $\bs{u} = \bs{v} + \bs{W}_{\! \! *}$ and $p = q + \Pi_{*}$ in $\Omega$, \eqref{stimacatta00} then follows directly from \eqref{fluxestimate}-\eqref{stima000}. Afterwards, a standard bootstrap argument enables us to deduce that $(\bs{u},p) \in W^{2,2}(\Omega;\mathbb{R}^N) \times W^{1,2}(\Omega;\mathbb{R})$.
\end{proof}

The main result of this article guarantees the unrestricted solvability (concerning the size of the data) of problem \eqref{nsstokes0}.
\begin{theorem} \label{epslevel}
Let $\Omega \subset \mathbb{R}^{N}$, $N \in \{2,3\}$, be an admissible domain. Given a force $\bs{f} \in L^2(\Omega;\mathbb{R}^N)$, $\bs{g}_{*} \in W^{3/2,2}_{+}(\Gamma_{I};\mathbb{R}^N)$ and $\sigma_{*} \in W^{1/2,2}(\Gamma_{O};\mathbb{R})$, there exists at least one weak solution $\bs{u} \in \mathcal{V}(\Omega)$ of problem \eqref{nsstokes0}, in the sense of Definition \ref{weaksolution}.
\end{theorem}
\noindent
\begin{proof}	
To prove the existence of a weak solution $\bs{u} \in \mathcal{V}(\Omega)$ of \eqref{nsstokes0} amounts to show the existence of $\widehat{\bs{u}} \in \mathcal{V}_{*}(\Omega)$ such that
\begin{equation}\label{oseenquasi}
\begin{aligned}
& \nu \int_{\Omega} \nabla \widehat{\bs{u}} : \nabla \bs{\varphi} + \int_{\Omega} (\widehat{\bs{u}} \cdot \nabla)\widehat{\bs{u}} \cdot \bs{\varphi} + \int_{\Omega} \left[ (\widehat{\bs{u}} \cdot \nabla)\bs{W}_{\! \! *} + (\bs{W}_{\! \! *} \cdot \nabla)\widehat{\bs{u}} \right] \cdot \bs{\varphi} + \dfrac{1}{2} \int_{\Gamma_{O}} [(\widehat{\bs{u}} + \bs{W}_{\! \! *}) \cdot \bs{n}]^{-}(\widehat{\bs{u}} \cdot \bs{\varphi}) \\[6pt]
& \hspace{-3mm}   = \int_{\Omega} \bs{f} \cdot \bs{\varphi} + \int_{\Gamma_{O}} \sigma_{*} (\bs{\varphi} \cdot \bs{n}) - \nu \int_{\Omega} \nabla \bs{W}_{\! \! *} : \nabla \bs{\varphi} - \int_{\Omega} (\bs{W}_{\! \! *} \cdot \nabla)\bs{W}_{\! \! *} \cdot \bs{\varphi} \qquad \forall \bs{\varphi} \in \mathcal{V}_{*}(\Omega) \, ,
\end{aligned}
\end{equation}
so that the solution will be given by $\bs{u}=\widehat{\bs{u}} + \bs{W}_{\! \! *}$. For a fixed $\widehat{\bs{u}} \in \mathcal{V}_{*}(\Omega)$, the applications
$$
\begin{aligned}
& \bs{\varphi} \in \mathcal{V}_{*}(\Omega) \longmapsto \int_{\Omega} (\widehat{\bs{u}} \cdot \nabla)\widehat{\bs{u}} \cdot \bs{\varphi} + \int_{\Omega} (\widehat{\bs{u}} \cdot \nabla)\bs{W}_{\! \! *} \cdot \bs{\varphi} + \int_{\Omega} (\bs{W}_{\! \! *} \cdot \nabla)\widehat{\bs{u}} \cdot \bs{\varphi} + \dfrac{1}{2} \int_{\Gamma_{O}} [(\widehat{\bs{u}} + \bs{W}_{\! \! *}) \cdot \bs{n}]^{-}(\widehat{\bs{u}} \cdot \bs{\varphi}) \, , \\[6pt]
& \bs{\varphi} \in \mathcal{V}_{*}(\Omega) \longmapsto \int_{\Omega} \bs{f} \cdot \bs{\varphi} + \int_{\Gamma_{O}} \sigma_{*} (\bs{\varphi} \cdot \bs{n}) - \nu \int_{\Omega} \nabla \bs{W}_{\! \! *} : \nabla \bs{\varphi} - \int_{\Omega} (\bs{W}_{\! \! *} \cdot \nabla)\bs{W}_{\! \! *} \cdot \bs{\varphi}
\end{aligned}
$$
clearly define linear continuous functionals on $\mathcal{V}_{*}(\Omega)$, see Remark \ref{rem1}. Then, in view of the Riesz Representation Theorem, the identity \eqref{oseenquasi} may be written as
$$
[\nu \, \widehat{\bs{u}} + \mathcal{P}(\widehat{\bs{u}}) - \mathcal{F}, \bs{\varphi}]_{\Omega} = 0 \qquad \forall \bs{\varphi} \in \mathcal{V}_{*}(\Omega) \, ,
$$
see \eqref{dirscal}, for some (unique) elements $\mathcal{P}(\widehat{\bs{u}})$, $\mathcal{F} \in \mathcal{V}_{*}(\Omega)$ such that
$$
\begin{aligned}
& [\mathcal{P}(\widehat{\bs{u}}), \bs{\varphi}]_{\Omega} = \int_{\Omega} (\widehat{\bs{u}} \cdot \nabla)\widehat{\bs{u}} \cdot \bs{\varphi} + \int_{\Omega} (\widehat{\bs{u}} \cdot \nabla)\bs{W}_{\! \! *} \cdot \bs{\varphi} + \int_{\Omega} (\bs{W}_{\! \! *} \cdot \nabla)\widehat{\bs{u}} \cdot \bs{\varphi} + \dfrac{1}{2} \int_{\Gamma_{O}} [(\widehat{\bs{u}} + \bs{W}_{\! \! *}) \cdot \bs{n}]^{-}(\widehat{\bs{u}} \cdot \bs{\varphi}) \, , \\[6pt]
& [\mathcal{F}, \bs{\varphi}]_{\Omega} = \int_{\Omega} \bs{f} \cdot \bs{\varphi} + \int_{\Gamma_{O}} \sigma_{*} (\bs{\varphi} \cdot \bs{n}) - \nu \int_{\Omega} \nabla \bs{W}_{\! \! *} : \nabla \bs{\varphi} - \int_{\Omega} (\bs{W}_{\! \! *} \cdot \nabla)\bs{W}_{\! \! *} \cdot \bs{\varphi} \qquad \forall \bs{\varphi} \in \mathcal{V}_{*}(\Omega) \, .
\end{aligned}
$$
We have so defined an operator $\mathcal{P} : \mathcal{V}_{*}(\Omega) \longrightarrow \mathcal{V}_{*}(\Omega)$ and we are led to find a solution $\widehat{\bs{u}} \in \mathcal{V}_{*}(\Omega)$ of the following nonlinear operator equation:
\begin{equation}\label{oseenquasi2}
	\widehat{\bs{u}} + \dfrac{1}{\nu}(\mathcal{P}(\widehat{\bs{u}}) - \mathcal{F}) = 0 \quad \text{in} \ \ \mathcal{V}_{*}(\Omega) \, .
\end{equation}
Exactly as in \cite[Chapter 5, Theorem 1]{ladyzhenskaya1969mathematical} one can show that the operator $\mathcal{P}$ is compact. Therefore, as a consequence of the Leray-Schauder Principle \cite[Chapter 6]{zeidler2013nonlinear}, in order to prove that \eqref{oseenquasi2} possesses at least one solution, it suffices to guarantee that any $\bs{v}^{\lambda} \in \mathcal{V}_{*}(\Omega)$ such that
\begin{equation}\label{oseenquasi3}
	\bs{v}^{\lambda} + \dfrac{\lambda}{\nu}(\mathcal{P}(\bs{v}^{\lambda}) - \mathcal{F}) = 0 \quad \text{in} \ \ \mathcal{V}_{*}(\Omega) \, ,
\end{equation}
is uniformly bounded with respect to $\lambda \in [0, 1]$. From now on, $C > 0$ will denote a generic constant that depends on $\Omega$ and $\nu$ (independent of $\lambda \in [0,1])$, but that may change from line to line.
\par
Given $\lambda \in [0, 1]$ and $\bs{v}^{\lambda} \in \mathcal{V}_{*}(\Omega) \setminus \{0\}$ such that \eqref{oseenquasi3} holds, we clearly have
\begin{equation} \label{vlambda0}
\begin{aligned}
& \nu \int_{\Omega} \nabla \bs{v}^{\lambda} : \nabla \bs{\varphi} + \lambda \int_{\Omega} \left[ (\bs{v}^{\lambda} \cdot \nabla)\bs{v}^{\lambda} + (\bs{v}^{\lambda} \cdot \nabla)\bs{W}_{\! \! *} + (\bs{W}_{\! \! *} \cdot \nabla)\bs{v}^{\lambda} \right] \cdot \bs{\varphi}  \\[6pt]
& + \dfrac{\lambda}{2} \int_{\Gamma_{O}} [(\bs{v}^{\lambda} + \bs{W}_{\! \! *}) \cdot \bs{n}]^{-}(\bs{v}^{\lambda} \cdot \bs{\varphi})  \\[6pt]
& \hspace{-4mm}  = \lambda \int_{\Omega} \bs{f} \cdot \bs{\varphi} + \lambda \int_{\Gamma_{O}} \sigma_{*} (\bs{\varphi} \cdot \bs{n}) - \lambda \nu \int_{\Omega} \nabla \bs{W}_{\! \! *} : \nabla \bs{\varphi} - \lambda \int_{\Omega} (\bs{W}_{\! \! *} \cdot \nabla)\bs{W}_{\! \! *} \cdot \bs{\varphi} \qquad \forall \bs{\varphi} \in \mathcal{V}_{*}(\Omega) \, .
\end{aligned}
\end{equation}
A standard integration by parts yields
\begin{equation} \label{vlambda1}
\begin{aligned}
\int_{\Omega} (\bs{v}^{\lambda} \cdot \nabla)\bs{v}^{\lambda} \cdot \bs{v}^{\lambda} = \dfrac{1}{2} \int_{\Gamma_{O}} | \bs{v}^{\lambda} |^{2}(\bs{v}^{\lambda} \cdot \bs{n}) \quad \text{and} \quad \int_{\Omega} (\bs{W}_{\! \! *} \cdot \nabla)\bs{v}^{\lambda} \cdot \bs{v}^{\lambda} = \dfrac{1}{2} \int_{\Gamma_{O}} | \bs{v}^{\lambda} |^{2}(\bs{W}_{\! \! *} \cdot \bs{n}) \, .
\end{aligned}
\end{equation}
After putting $\bs{\varphi} = \bs{v}^{\lambda}$ in \eqref{vlambda0} and enforcing \eqref{vlambda1}, the following identity is obtained:
$$
\begin{aligned}
\nu \| \nabla \bs{v}^{\lambda} \|^{2}_{L^{2}(\Omega)} & = \lambda \int_{\Omega} \bs{f} \cdot \bs{v}^{\lambda} + \lambda \int_{\Gamma_{O}} \sigma_{*} (\bs{v}^{\lambda} \cdot \bs{n}) - \lambda \nu \int_{\Omega} \nabla \bs{W}_{\! \! *} : \nabla \bs{v}^{\lambda} - \lambda \int_{\Omega} (\bs{W}_{\! \! *} \cdot \nabla)\bs{W}_{\! \! *} \cdot \bs{v}^{\lambda}  \\[6pt]
& \hspace{4mm} - \lambda \int_{\Omega} (\bs{v}^{\lambda} \cdot \nabla)\bs{W}_{\! \! *} \cdot \bs{v}^{\lambda} - \dfrac{\lambda}{2} \int_{\Gamma_{O}} \left( (\bs{v}^{\lambda} + \bs{W}_{\! \! *}) \cdot \bs{n} + [(\bs{v}^{\lambda} + \bs{W}_{\! \! *}) \cdot \bs{n}]^{-} \right) | \bs{v}^{\lambda} |^{2} \, ,
\end{aligned}
$$
so that
\begin{equation} \label{important}
\begin{aligned}	
\nu \| \nabla \bs{v}^{\lambda} \|^{2}_{L^{2}(\Omega)} & \leq \lambda \int_{\Omega} \bs{f} \cdot \bs{v}^{\lambda} + \lambda \int_{\Gamma_{O}} \sigma_{*} (\bs{v}^{\lambda} \cdot \bs{n}) - \lambda \nu \int_{\Omega} \nabla \bs{W}_{\! \! *} : \nabla \bs{v}^{\lambda} - \lambda \int_{\Omega} (\bs{W}_{\! \! *} \cdot \nabla)\bs{W}_{\! \! *} \cdot \bs{v}^{\lambda}  \\[6pt]
& \hspace{4mm} - \lambda \int_{\Omega} (\bs{v}^{\lambda} \cdot \nabla)\bs{W}_{\! \! *} \cdot \bs{v}^{\lambda} \qquad \forall \lambda \in [0,1] \, . 
\end{aligned} 
\end{equation}
The integral identity \eqref{vlambda0} and Theorem \ref{press} guarantee the existence of a unique pressure $p^{\lambda} \in L^{2}(\Omega; \mathbb{R})$ such that 
\begin{equation}\label{nsstokeslambdapre}
	\begin{aligned}
		& \nu \int_{\Omega} \nabla \bs{v}^{\lambda} : \nabla \bs{\varphi} + \lambda \int_{\Omega} \left[ (\bs{v}^{\lambda} \cdot \nabla)\bs{v}^{\lambda} + (\bs{v}^{\lambda} \cdot \nabla)\bs{W}_{\! \! *} + (\bs{W}_{\! \! *} \cdot \nabla)\bs{v}^{\lambda} \right] \cdot \bs{\varphi} - \int_{\Omega} p^{\lambda} (\nabla \cdot \bs{\varphi})  \\[6pt]
		&  + \dfrac{\lambda}{2} \int_{\Gamma_{O}} [(\bs{v}^{\lambda} + \bs{W}_{\! \! *}) \cdot \bs{n}]^{-}(\bs{v}^{\lambda} \cdot \bs{\varphi})  \\[6pt]
		& \hspace{-4mm}  =  \lambda \int_{\Omega} \bs{f} \cdot \bs{\varphi} + \lambda \int_{\Gamma_{O}} \sigma_{*} (\bs{\varphi} \cdot \bs{n}) - \lambda \nu \int_{\Omega} \nabla \bs{W}_{\! \! *} : \nabla \bs{\varphi} - \lambda \int_{\Omega} (\bs{W}_{\! \! *} \cdot \nabla)\bs{W}_{\! \! *} \cdot \bs{\varphi} \qquad \forall \bs{\varphi} \in \mathcal{H}_{*}(\Omega) \, .
	\end{aligned}	
\end{equation}
Define $q^{\lambda} \doteq p^{\lambda} - \lambda \, \Pi_{*}$ in $\Omega$, for every $\lambda \in [0,1]$, so that \eqref{nsstokeslambdapre} becomes
\begin{equation}\label{nsstokeslambda}
\begin{aligned}
	& \nu \int_{\Omega} \nabla \bs{v}^{\lambda} : \nabla \bs{\varphi} + \lambda \int_{\Omega} \left[ (\bs{v}^{\lambda} \cdot \nabla)\bs{v}^{\lambda} + (\bs{v}^{\lambda} \cdot \nabla)\bs{W}_{\! \! *} + (\bs{W}_{\! \! *} \cdot \nabla)\bs{v}^{\lambda} \right] \cdot \bs{\varphi} - \int_{\Omega} q^{\lambda} (\nabla \cdot \bs{\varphi})  \\[6pt]
	&  + \dfrac{\lambda}{2} \int_{\Gamma_{O}} [(\bs{v}^{\lambda} + \bs{W}_{\! \! *}) \cdot \bs{n}]^{-}(\bs{v}^{\lambda} \cdot \bs{\varphi})  \\[6pt]
	& \hspace{-4mm}  =  \lambda \int_{\Omega} \bs{f} \cdot \bs{\varphi} + \lambda \int_{\Gamma_{O}} \sigma_{*} (\bs{\varphi} \cdot \bs{n}) - \lambda \nu \int_{\Omega} \nabla \bs{W}_{\! \! *} : \nabla \bs{\varphi} - \lambda \int_{\Omega} (\bs{W}_{\! \! *} \cdot \nabla)\bs{W}_{\! \! *} \cdot \bs{\varphi} + \lambda \int_{\Omega} \Pi_{*} (\nabla \cdot \bs{\varphi}) \, ,
\end{aligned}	
\end{equation}
for every $\bs{\varphi} \in \mathcal{H}_{*}(\Omega)$. Moreover, $(\bs{v}^{\lambda},q^{\lambda}) \in W^{2,2}(\Omega;\mathbb{R}^N) \times W^{1,2}(\Omega;\mathbb{R})$, so that the system
\begin{equation}\label{nsstokeslambdaforte}
	\left\{
	\begin{aligned}
		& -\nu\Delta \bs{v}^{\lambda} + \lambda \left[ (\bs{v}^{\lambda} \cdot \nabla)\bs{v}^{\lambda} + (\bs{v}^{\lambda} \cdot \nabla)\bs{W}_{\! \! *} + (\bs{W}_{\! \! *} \cdot \nabla)\bs{v}^{\lambda} \right] + \nabla q^{\lambda} \\[5pt]
		& \hspace{-2mm} = \lambda \left[ \bs{f} + \nu\Delta \bs{W}_{\! \! *} - (\bs{W}_{\! \! *} \cdot \nabla)\bs{W}_{\! \! *}  - \nabla \Pi_{*} \right] \quad \text{in} \ \ \Omega \, , \\[5pt]
& \nabla \cdot \bs{v}^{\lambda} = 0 \quad \text{in} \ \ \Omega \, , \\[5pt]
& \bs{v}^{\lambda} = \bs{0} \quad \text{on} \ \ \Gamma_{I} \cup \Gamma_{W} \, , \\[5pt]
		& \nu \dfrac{\partial \bs{v}^{\lambda}}{\partial \bs{n}} - q^{\lambda} \bs{n} + \dfrac{\lambda}{2} [(\bs{v}^{\lambda} + \bs{W}_{\! \! *}) \cdot \bs{n}]^{-}\bs{v}^{\lambda} = \bs{0} \quad \mbox{on} \ \ \Gamma_{O} \, ,
	\end{aligned}
	\right.
\end{equation}
is satisfied in strong form, alongside the following estimate (see \eqref{stima000}):
\begin{equation} \label{stokesest2}
\begin{aligned}
\| \bs{v}^{\lambda} \|_{W^{2,3/2}(\Omega)} + \| q^{\lambda} \|_{W^{1,3/2}(\Omega)} & \leq C \Big( \| \bs{f} \|_{L^{2}(\Omega)} + \| \nabla \bs{v}^{\lambda} \|^{2}_{L^{2}(\Omega)} + \| \nabla \bs{W}_{\! \! *} \|^{2}_{L^{2}(\Omega)}  \\[6pt] 
	& \hspace{1.2cm} + \| \bs{W}_{\! \! *} \|_{W^{2,2}(\Omega)} + \| \Pi_{*} \|_{W^{1,2}(\Omega)} \Big) \qquad \forall \lambda \in [0,1] \, .
\end{aligned}
\end{equation}
By contradiction, suppose now that the norms $\| \nabla \bs{v}^{\lambda} \|_{L^{2}(\Omega)}$ are \textit{not} uniformly bounded with respect to $\lambda \in [0,1]$. Then, there must exist $\lambda_{0} \in [0,1]$ and a sequence $(\lambda_{k})_{k \in \mathbb{N}} \subset [0,1]$ such that
\begin{equation} \label{contra1}
\lim_{k \to \infty} \lambda_{k} = \lambda_{0} \qquad \text{and} \qquad \lim_{k \to \infty} J_{k} = + \infty \, , \quad \text{with} \ \ J_{k} \doteq \| \nabla \bs{v}^{\lambda_{k}} \|_{L^{2}(\Omega)} \quad \forall k \in \mathbb{N} \, .
\end{equation}
The estimate in \eqref{stokesest2} entails that, along the divergent sequence \eqref{contra1}, the following sequences are all uniformly bounded with respect to $k \in \mathbb{N}$:
$$
\begin{aligned}
& ( \widehat{\bs{v}}_{k} )_{k \in \mathbb{N}} \doteq \left( \dfrac{\bs{v}^{\lambda_{k}}}{J_{k}} \right)_{k \in \mathbb{N}} \subset \mathcal{V}_{*}(\Omega) \, ; \qquad ( \widetilde{\bs{v}}_{k} )_{k \in \mathbb{N}} \doteq \left( \dfrac{\bs{v}^{\lambda_{k}}}{J^{2}_{k}} \right)_{k \in \mathbb{N}} \subset W^{2,3/2}(\Omega;\mathbb{R}^N) \cap \mathcal{V}_{*}(\Omega) \, ; \\[3pt]
& \hspace{3cm} ( \widehat{q}_{k} )_{k \in \mathbb{N}} \doteq  \left( \dfrac{q^{\lambda_k}}{J^{2}_{k}} \right)_{k \in \mathbb{N}} \subset W^{1,3/2}(\Omega;\mathbb{R}) \, .
\end{aligned}
$$
Therefore, there exist $\widehat{\bs{v}} \in \mathcal{V}_{*}(\Omega)$, $\widetilde{\bs{v}} \in W^{2,3/2}(\Omega;\mathbb{R}^N) \cap \mathcal{V}_{*}(\Omega)$ and $\widehat{q} \in W^{1,3/2}(\Omega;\mathbb{R})$ such that the following convergences hold as $k \to \infty$:
\begin{equation} \label{convergences1}
\begin{aligned}
& \widehat{\bs{v}}_{k} \rightharpoonup \widehat{\bs{v}} \ \ \ \text{weakly in} \ \ \ \mathcal{V}_{*}(\Omega) \, ; \qquad \widehat{\bs{v}}_{k} \to \widehat{\bs{v}} \ \ \ \text{in} \ \ \ L^{4}(\Omega;\mathbb{R}^N) \, ; \qquad \widehat{\bs{v}}_{k} \to \widehat{\bs{v}} \ \ \ \text{in} \ \ \ L^{2}(\partial \Omega;\mathbb{R}^N) \, ; \\[3pt]
& \widetilde{\bs{v}}_{k} \to \widetilde{\bs{v}} \ \ \ \text{in} \ \ \ W^{1,3/2}(\Omega;\mathbb{R}^N) \, ; \qquad \dfrac{\partial \widetilde{\bs{v}}_{k}}{\partial \bs{n}} \to \dfrac{\partial \widetilde{\bs{v}}}{\partial \bs{n}} \ \ \ \text{in} \ \ \ L^{1}(\Gamma_{O};\mathbb{R}^N) \, ; \\[3pt]
& \widehat{q}_{k} \rightharpoonup \widehat{q} \ \ \ \text{weakly in} \ \ \ W^{1,3/2}(\Omega;\mathbb{R}) \, ; \qquad \widehat{q}_{k} \to \widehat{q} \ \ \ \text{in} \ \ \ L^{1}(\partial \Omega;\mathbb{R}) \, ; \\[3pt]
\end{aligned}
\end{equation}
along subsequences that are not being relabeled, see also \cite[Theorem 6.2]{necas2011direct}. An application of the H\"older inequality readily implies that
\begin{equation} \label{convergences11} 
\widetilde{\bs{v}} \equiv \bs{0} \ \ \ \text{in} \ \ \ \Omega \, ,
\end{equation}
since $\| \nabla \widetilde{\bs{v}}_{k} \|_{L^{2}(\Omega)} \to 0$ as $k \to \infty$. Given $k \in \mathbb{N}$, if we divide the inequality \eqref{important} (with $\lambda = \lambda_{k}$) by $J^{2}_{k}$, we obtain
\begin{equation} \label{vlambda4}
\begin{aligned}
\nu  & \leq \dfrac{\lambda_{k}}{J_{k}} \int_{\Omega} \bs{f} \cdot \widehat{\bs{v}}_{k} + \dfrac{\lambda_{k}}{J_{k}} \int_{\Gamma_{O}} \sigma_{*} (\widehat{\bs{v}}_{k} \cdot \bs{n}) - \nu \dfrac{\lambda_{k}}{J_{k}} \int_{\Omega} \nabla \bs{W}_{\! \! *} : \nabla \widehat{\bs{v}}_{k} - \dfrac{\lambda_{k}}{J_{k}} \int_{\Omega} (\bs{W}_{\! \! *} \cdot \nabla)\bs{W}_{\! \! *} \cdot \widehat{\bs{v}}_{k}  \\[6pt]
& \hspace{4mm} - \lambda_{k} \int_{\Omega} (\widehat{\bs{v}}_{k} \cdot \nabla)\bs{W}_{\! \! *} \cdot \widehat{\bs{v}}_{k} \qquad \forall k \in \mathbb{N} \, ,
\end{aligned}
\end{equation}
where we may then take the limit as $k \to +\infty$ (invoking \eqref{convergences1}) and deduce
\begin{equation} \label{vlambdalim}
\nu \leq - \lambda_{0} \int_{\Omega} (\widehat{\bs{v}} \cdot \nabla)\bs{W}_{\! \! *} \cdot \widehat{\bs{v}} \, ,
\end{equation}
thus implying that $\lambda_{0} > 0$. On the other hand, given any $k \in \mathbb{N}$ and $\bs{\varphi} \in \mathcal{C}^{\infty}_{0}(\Omega;\mathbb{R}^N) \subset \mathcal{H}_{*}(\Omega)$ (not necessarily divergence-free), after testing identity \eqref{nsstokeslambda} (with $\lambda = \lambda_{k}$) with $\bs{\varphi}$ (so that all boundary integrals vanish), and dividing the resulting equality by $J^{2}_{k}$, we obtain
\begin{equation} \label{vlambdalim2}
\begin{aligned}
& \dfrac{\nu}{J_{k}} \int_{\Omega} \nabla \widehat{\bs{v}}_{k} : \nabla \bs{\varphi} + \lambda_{k} \int_{\Omega} (\widehat{\bs{v}}_{k} \cdot \nabla)\widehat{\bs{v}}_{k} \cdot \bs{\varphi} + \dfrac{\lambda_{k}}{J_{k}}  \int_{\Omega} (\widehat{\bs{v}}_{k} \cdot \nabla)\bs{W}_{\! \! *} \cdot \bs{\varphi} + \dfrac{\lambda_{k}}{J_{k}}  \int_{\Omega} (\bs{W}_{\! \! *} \cdot \nabla)\widehat{\bs{v}}_{k} \cdot \bs{\varphi} \\[6pt]
& \hspace{-4mm} = \dfrac{\lambda_{k}}{J^{2}_{k}} \int_{\Omega} \bs{f} \cdot \bs{\varphi} - \nu \dfrac{\lambda_{k}}{J^{2}_{k}} \int_{\Omega} \nabla \bs{W}_{\! \! *} : \nabla \bs{\varphi} - \dfrac{\lambda_{k}}{J^{2}_{k}} \int_{\Omega} (\bs{W}_{\! \! *} \cdot \nabla)\bs{W}_{\! \! *} \cdot \bs{\varphi} + \int_{\Omega} \left( \widehat{q}_{k} + \lambda_{k} \dfrac{\Pi_{*}}{J^{2}_{k}} \right) (\nabla \cdot \bs{\varphi}) \, ,
\end{aligned}
\end{equation}
for every $k \in \mathbb{N}$, where we may then take the limit as $k \to +\infty$ (invoking \eqref{convergences1}) and deduce that
\begin{equation} \label{limiteulerdis}
\lambda_{0} \int_{\Omega} (\widehat{\bs{v}} \cdot \nabla)\widehat{\bs{v}} \cdot \bs{\varphi} - \int_{\Omega} \widehat{q} (\nabla \cdot \bs{\varphi})  = 0 \qquad \forall \bs{\varphi} \in \mathcal{C}^{\infty}_{0}(\Omega;\mathbb{R}^N) \, ,
\end{equation}
that is, the pair $(\widehat{\bs{v}}, \widehat{q}) \in \mathcal{V}_{*}(\Omega) \times W^{1,3/2}(\Omega;\mathbb{R})$ satisfies in strong form the following Euler-type equation:
\begin{equation} \label{limiteuler}
\lambda_{0} (\widehat{\bs{v}} \cdot \nabla)\widehat{\bs{v}} + \nabla \widehat{q} = \bs{0} \, , \quad  \nabla\cdot \widehat{\bs{v}}=0 \ \ \mbox{ in } \ \ \Omega \, .
\end{equation}
By setting $\widehat{\bs{v}}_{0} \doteq \sqrt{\lambda_{0}} \, \widehat{\bs{v}}$, we notice that the pair $(\widehat{\bs{v}}_{0}, \widehat{q}) \in \mathcal{V}_{*}(\Omega) \times W^{1,3/2}(\Omega;\mathbb{R})$ is a strong solution of the Euler equation
$$
(\widehat{\bs{v}}_{0} \cdot \nabla) \widehat{\bs{v}}_{0} + \nabla \widehat{q} = \bs{0} \, , \quad  \nabla\cdot \widehat{\bs{v}}_{0}=0 \ \ \mbox{ in } \ \ \Omega \, .
$$
Since $\widehat{\bs{v}}_{0} = \bs{0}$ on $\Gamma_{I} \cup \Gamma_{W}$ (which is one connected component of $\partial \Omega)$, the Bernoulli law \cite[Lemma 4]{kapitanskii1983spaces} (see \cite[Theorem 2.2]{amick1984existence} and \cite[Theorem 1]{korobkov2011bernoulli} as well) states that $\widehat{q}$ must be constant on $\Gamma_{I} \cup \Gamma_{W}$: there exists $\widehat{q}_{*} \in \mathbb{R}$ such that
\begin{equation} \label{bernoulli1}
\widehat{q} = \widehat{q}_{*} \quad \text{almost everywhere on} \ \ \Gamma_{I} \cup \Gamma_{W} \, .
\end{equation}
On the other hand, notice that \eqref{nsstokeslambdaforte}$_4$ ensures
$$
\widehat{q}_{k} \, \bs{n} = \nu \dfrac{\partial \widetilde{\bs{v}}_{k}}{\partial \bs{n}}  + \dfrac{\lambda_{k}}{2} \left[ \left( \widehat{\bs{v}}_{k} + \dfrac{\bs{W}_{\! \! *}}{J_{k}} \right) \cdot \bs{n} \right]^{-} \widehat{\bs{v}}_{k}  \quad \mbox{a.e. on} \ \ \Gamma_{O} \, , \quad \forall k \in \mathbb{N} \, ,
$$
where we may then take the limit as $k \to +\infty$ (invoking \eqref{convergences1}-\eqref{convergences11}) and conclude that
\begin{equation} \label{bernoulli2}
	  \widehat{q} \, \bs{n} = \dfrac{\lambda_{0}}{2} [\widehat{\bs{v}} \cdot \bs{n}]^{-} \, \widehat{\bs{v}} \quad \mbox{a.e. on} \ \ \Gamma_{O}  \, ,
\end{equation}
that is,
\begin{equation} \label{bernoulli22}
\widehat{q}  = \dfrac{\lambda_{0}}{2} [\widehat{\bs{v}} \cdot \bs{n}]^{-} \, (\widehat{\bs{v}} \cdot \bs{n}) \leq 0 \quad \mbox{a.e. on} \ \ \Gamma_{O} \, .
\end{equation}

We then claim that $\widehat{q}_{*}\leq 0$; to prove this, we adapt the arguments of \cite[Theorem 3.1]{korobkov2020solvability} to our current configuration in the case when $N=3$ (namely, cylindrical ends with general cross sections), since the case $N=2$ was already covered in \cite[Theorem 3.1]{korobkov2020solvability}. Denote by $\bs{n}(\bs{s}) \in \mathbb{R}^{2}$ the outward unit normal to $\partial \Theta_{2}$ at a point $\bs{s} \in \partial \Theta_{2}$, so that $\bs{n} \in \mathcal{C}^{1}(\partial \Theta_{2}; \mathbb{R}^{2})$, as $\partial \Theta_{2}$ is of class $\mathcal{C}^{2}$. There exists a sufficiently small $\eps_{*} > 0$ (say $0< \eps_{*} \ll \text{diam}(\Theta_{2}))$ such that, after putting
\begin{equation} \label{bl0}
\Theta^{\eps}_{2} \doteq \{\bs{y}\in\Theta_2 \ | \ \text{dist}(\bs{y},\partial\Theta_{2}) < \eps \, \} \qquad \forall \eps \in (0,\eps_{*}] \, ,
\end{equation} 
the mapping $\Psi : \partial \Theta_{2} \times [0,\eps_{*}] \longrightarrow \overline{\Theta^{\eps_{*}}_{2}}$ defined by
\begin{equation} \label{bl1}
\Psi(\bs{s},t) = \bs{y} \doteq \bs{s} - t \bs{n}(\bs{s}) \qquad \forall (\bs{s},t) \in \partial\Theta_{2} \times [0,\eps_{*}] \, ,
\end{equation}
is a $\mathcal{C}^{1}$-diffeomorphism whose Jacobian can be uniformly bounded (also away from zero) with respect to $\eps \in (0, \eps_{*})$. Then, given $\eps \in (0, \eps_{*})$, define the \textit{boundary layer} of width $\eps$, included in $\Omega_2$, as
\begin{equation} \label{bl3}
\Omega^{\eps}_{2} \doteq \Theta^{\eps}_{2} \times (0,\ell_2) = \Psi(\partial \Theta_{2} \times (0,\eps))  \times (0,\ell_2) \, ,
\end{equation}
as well as the parallel surfaces
\begin{equation} \label{bl4}
S_{\eps} \doteq \{\bs{y} \in \Theta_2 \ | \ \text{dist}(\bs{y},\partial\Theta_{2}) = \eps \, \} = \Psi(\partial \Theta_{2}\times \{ \eps \}) \, .
\end{equation}
Since $\widehat{q} \in W^{1,3/2}(\Omega;\mathbb{R})$, for almost every $(\bs{s},\eps,z) \in \partial\Theta_{2} \times (0,\eps_{*}) \times (0,\ell_{2})$ we may write
$$
\begin{aligned}
\widehat{q}(\Psi(\bs{s},\eps),z) - \widehat{q}_{*} & = \widehat{q}(\Psi(\bs{s},\eps),z) - \widehat{q}(\Psi(\bs{s},0),z) \\[6pt]
& = \int_{0}^{\eps} \dfrac{\partial}{\partial \tau} \widehat{q}(\Psi(\bs{s},\tau),z) \, d\tau \\[6pt]
& = - \int_{0}^{\eps} \nabla_{\! \bs{y}} \widehat{q}(\Psi(\bs{s},\tau),z) \cdot \bs{n}(\bs{s}) \, d\tau \, ,
\end{aligned}
$$
thereby implying
$$
| \widehat{q}(\Psi(\bs{s},\eps),z) - \widehat{q}_{*} | \leq \int_{0}^{\eps} | \nabla \widehat{q}(\Psi(\bs{s},\tau),z) | \, d\tau \qquad \forall (\bs{s},\eps,z) \in \partial\Theta_{2} \times (0,\eps_{*}) \times (0,\ell_{2})-\text{a.e.}.
$$
Upon integration of this last inequality over $(\bs{s},z) \in \partial\Theta_{2} \times (0,\ell_{2})$ we notice, as a consequence of \eqref{bl3}-\eqref{bl4} and a change of variables, that
\begin{equation}
\int\limits_{S_{\eps} \times (0,\ell_2)} | \widehat{q} - \widehat{q}_{*} | \leq C \int_{\Omega^{\eps}_{2}}|\nabla\widehat{q}| \qquad \forall \eps \in (0, \eps_{*}) \, ,
\end{equation}
where $C>0$ is independent of $\eps \in (0, \eps_{*})$. It then follows
\begin{equation} \label{firststep}
\lim_{\eps \to 0^{+}} \int\limits_{S_{\eps} \times (0,\ell_2)} | \widehat{q} - \widehat{q}_{*} | = 0 \, .
\end{equation}
On the other hand, recalling that $\widehat{\bs{v}}$ vanishes on $\partial\Theta_2\times (0,\ell_2)$, the following Poincaré inequality holds:
\begin{equation}\label{poincare}
	\|\widehat{\bs{v}}\|_{L^2(\Omega^{\eps}_{2})}\leq C \eps \|\nabla\widehat{\bs{v}}\|_{L^2(\Omega^{\eps}_{2})} \qquad \forall \eps \in (0, \eps_{*}) \, ,
\end{equation}
where $C>0$ is independent of $\eps \in (0, \eps_{*})$. As the Euler equation \eqref{limiteuler} holds in $\Omega^{\eps}_{2}$, from the H\"older inequality and \eqref{poincare} we further infer
\begin{equation} \label{eulereps}
	\int_{\Omega^{\eps}_{2}}|\nabla\widehat{q}|=\lambda_0 \int_{\Omega^{\eps}_{2}}|(\widehat{\bs{v}}\cdot\nabla)\widehat{\bs{v}}|\leq \lambda_0\|\widehat{\bs{v}}\|_{L^2(\Omega^{\eps}_{2})}\|\nabla\widehat{\bs{v}}\|_{L^2(\Omega^{\eps}_{2})}\leq C\eps\|\nabla\widehat{\bs{v}}\|_{L^2(\Omega^{\eps}_{2})}^2 \qquad \forall \eps \in (0, \eps_{*}) \, .
\end{equation}
Now, the \textit{Coarea formula} \cite[Theorem 1.5.17]{korobkov2024steady} entails
$$
\int\limits_{0}^{\eps} \left( \, \int\limits_{S_{t} \times (0,\ell_2)} | \nabla\widehat{q} | \right) dt = \int_{\Omega^{\eps}_{2}}|\nabla\widehat{q}| \qquad \forall \eps \in (0, \eps_{*}) \, ,
$$
and thus, the Mean Value Theorem for Lebesgue Integrals ensures
$$
\left| \left\{ t \in (0, \eps) \ \Bigg| \ \int\limits_{S_{t} \times (0,\ell_2)} | \nabla\widehat{q} | \leq \dfrac{1}{\eps} \int_{\Omega^{\eps}_{2}}|\nabla\widehat{q}| \, \right\} \right| > 0 \qquad \forall \eps \in (0, \eps_{*}) \, .
$$
Therefore, for every $\eps \in (0, \eps_{*})$ we can find $\delta_{\eps} \in (0,\eps)$ such that
\begin{equation} \label{deltaeps1}
\int\limits_{S_{\delta_{\eps}} \times (0,\ell_2)} | \nabla\widehat{q} | \leq \dfrac{1}{\eps} \int_{\Omega^{\eps}_{2}}|\nabla\widehat{q}| \leq C \|\nabla\widehat{\bs{v}}\|_{L^2(\Omega^{\eps}_{2})}^2 \, .
\end{equation}
as a consequence of \eqref{eulereps}. Now, take $\eps \in (0, \eps_{*})$ and $\delta_{\eps} \in (0,\eps)$ satisfying \eqref{deltaeps1}. For almost every $(\bs{y},z) \in S_{\delta_{\eps}} \times (0,\ell_2)$ we may integrate in local coordinates to obtain
\begin{equation} \label{deltaeps2}
\widehat{q}(\bs{y},z) = \widehat{q}(\bs{y},\ell_{2}) - \int_{z}^{\ell_{2}} \dfrac{\partial \widehat{q}}{\partial z}(\bs{y},\tau) \, d\tau \, .
\end{equation}
Since $(\bs{y},\ell_{2}) \in \Gamma_{O}$, from \eqref{bernoulli22}-\eqref{deltaeps2} we deduce
$$
\widehat{q}(\bs{y},z) \leq - \int_{z}^{\ell_{2}} \dfrac{\partial \widehat{q}}{\partial z}(\bs{y},\tau) \, d\tau \leq \int_{0}^{\ell_{2}} | \nabla \widehat{q}(\bs{y},\tau) | d\tau \qquad \forall (\bs{y},z) \in S_{\delta_{\eps}} \times (0,\ell_2)-\text{a.e.} \, .
$$
Upon integration of this last inequality over $S_{\delta_{\eps}} \times (0,\ell_2)$ we notice, in view of \eqref{deltaeps1}, that
\begin{equation} \label{finalp0}
\int\limits_{S_{\delta_{\eps}} \times (0,\ell_2)} \widehat{q} \leq \ell_{2} \int\limits_{S_{\delta_{\eps}} \times (0,\ell_2)} | \nabla\widehat{q} | \leq \dfrac{\ell_2}{\eps} \int_{\Omega^{\eps}_{2}}|\nabla\widehat{q}| \leq C \|\nabla\widehat{\bs{v}}\|_{L^2(\Omega^{\eps}_{2})}^2 \qquad \forall \eps \in (0, \eps_{*})  \, .
\end{equation} 
In the light of \eqref{firststep}, we can take the limit as $\eps \to 0^{+}$ on both sides of inequality \eqref{finalp0} to deduce
$$
\ell_{2} | \partial \Theta_{2} |  \widehat{q}_{*} \leq 0
$$
or equivalently,
\begin{equation} \label{finalp}
	\widehat{q}_{*} \leq 0 \, .
\end{equation} 

Recall from Lemma \ref{lemma0} that $\bs{W}_{\! \! *} = \bs{V}^{(2)}_{\! \Phi_{*}}$ on $\Gamma_{O}$, see again \eqref{poi23dd}. In view of \eqref{poi1}-\eqref{maxprin}-\eqref{poi23d}-\eqref{poi23dd}, we deduce the existence of a non-negative function $\psi_{*} \in \mathcal{C}(\Gamma_{O};[0,+\infty))$ such that
\begin{equation} \label{finalp2}
\bs{W}_{\! \! *} = \psi_{*} \, \bs{n} \quad \mbox{on} \ \ \Gamma_{O} \, .
\end{equation} 
Multiply \eqref{limiteuler} by $\bs{W}_{\! \! *}$ and integrate by parts in $\Omega$ (enforcing \eqref{bernoulli1}-\eqref{bernoulli22}-\eqref{finalp}-\eqref{finalp2}, and recalling that $2 z^2 + [z]^{-} z \geq 0$, for every $z \in \mathbb{R}$), to obtain
$$
\begin{aligned}
- \lambda_{0} \int_{\Omega} (\widehat{\bs{v}} \cdot \nabla)\bs{W}_{\! \! *} \cdot \widehat{\bs{v}} & = \widehat{q}_{*} \Phi_{*}  - \int_{\Gamma_{O}} \widehat{q} (\bs{W}_{\! \! *} \cdot \bs{n}) - \lambda_{0} \int_{\Gamma_{O}} (\widehat{\bs{v}} \cdot \bs{W}_{\! \! *})(\widehat{\bs{v}} \cdot \bs{n}) \\[6pt]
& \leq -\lambda_{0} \int_{\Gamma_{O}}  \left( \dfrac{1}{2} [\widehat{\bs{v}} \cdot \bs{n}]^{-} \, (\widehat{\bs{v}} \cdot \bs{n}) + (\widehat{\bs{v}} \cdot \bs{n})^{2} \right) \psi_{*} \, ,
\end{aligned}
$$
that is,
$$
- \lambda_{0} \int_{\Omega} (\widehat{\bs{v}} \cdot \nabla)\bs{W}_{\! \! *} \cdot \widehat{\bs{v}} \leq 0 \, ,
$$
which certainly contradicts \eqref{vlambdalim}. Therefore, the norms $\| \nabla \bs{v}^{\lambda} \|_{L^{2}(\Omega)}$ must be uniformly bounded with respect to $\lambda \in [0,1]$, and so there exists $\widehat{\bs{u}} \in \mathcal{V}_{*}(\Omega)$ satisfying \eqref{oseenquasi}. This concludes the proof.
\end{proof}

Theorem \ref{epslevel} deserves some remarks and comments.

\begin{remark} \label{2dpreszero}
In the case when the admissible domain $\Omega \subset \mathbb{R}^{2}$ is a horizontal rectangle, for example, $\Omega = (0,\ell) \times (-h,h)$ (for some $\ell > 0$ and $h > 0$), we can actually show that $\widehat{q}_{*} = 0$, see again \eqref{finalp}. In fact, writing $\widehat{\bs{v}} = (\widehat{v}_{1}, \widehat{v}_{2})$ in $\Omega$, the first component of the Euler-type equation \eqref{limiteuler} reads
$$
\dfrac{\partial \widehat{q}}{\partial x} = - \lambda_{0} \left( \widehat{v}_{1} \dfrac{\partial \widehat{v}_{1}}{\partial x} + \widehat{v}_{2} \dfrac{\partial \widehat{v}_{1}}{\partial y} \right) \quad \mbox{a.e. in} \ \ \Omega \, .
$$
Integrating this last identity in $\Omega$ (integrating by parts the term on the right-hand side) entails
$$
\int_{-h}^{h} \left( \widehat{q}(\ell,y) - \widehat{q}_{*} \right) dy = - \lambda_{0} \int_{-h}^{h} | \widehat{v}_{1}(\ell,y) |^{2} \, dy \, ,
$$ 
which, owing to \eqref{bernoulli22}, means that
$$
\widehat{q}_{*} = \dfrac{\lambda_{0}}{2h} \int_{-h}^{h} \left( \dfrac{1}{2} [\widehat{v}_{1}(\ell,y)]^{-} \, \widehat{v}_{1}(\ell,y) + | \widehat{v}_{1}(\ell,y) |^{2} \right) dy \geq 0 \, .
$$
Together with \eqref{finalp}, this last identity implies that $\widehat{q}_{*} = 0$.
\end{remark}

\begin{remark} \label{3dpreszero}
In analogy with Remark \ref{2dpreszero}, in the case when the admissible domain $\Omega \subset \mathbb{R}^{3}$ is a horizontal pipe, say, $\Omega = \Theta \times (0,\ell)$ (for some $\ell>0$ and some open, bounded and planar domain $\Theta \subset \mathbb{R}^{2}$ of class $\mathcal{C}^{2}$), we can also show that $\widehat{q}_{*} = 0$, see again \eqref{finalp}. In fact, writing $\widehat{\bs{v}} = (\widehat{v}_{1}, \widehat{v}_{2}, \widehat{v}_{3} )$ in $\Omega$, the third component of the Euler-type equation \eqref{limiteuler} reads
	$$
	\dfrac{\partial \widehat{q}}{\partial z} = - \lambda_{0} \left( \widehat{v}_{1} \dfrac{\partial \widehat{v}_{3}}{\partial x} + \widehat{v}_{2} \dfrac{\partial \widehat{v}_{3}}{\partial y} + \widehat{v}_{3} \dfrac{\partial \widehat{v}_{3}}{\partial z} \right) \quad \mbox{a.e. in} \ \ \Omega \, .
	$$
	Integrating this identity in $\Omega$ (integrating by parts the term on the right-hand side) entails
	$$
	\int_{\Theta} \left( \widehat{q}(x,y,\ell) - \widehat{q}_{*} \right) dx \, dy = - \lambda_{0} \int_{\Theta} | \widehat{v}_{3}(x,y,\ell) |^{2} \, dx \, dy \, ,
	$$ 
	which, owing to \eqref{bernoulli22}, means that
	$$
	\widehat{q}_{*} = \dfrac{\lambda_{0}}{| \Theta |} \int_{\Theta} \left( \dfrac{1}{2} [\widehat{v}_{3}(x,y,\ell)]^{-} \, \widehat{v}_{3}(x,y,\ell) + | \widehat{v}_{3}(x,y,\ell) |^{2} \right) dx \, dy \geq 0 \, .
	$$
	Together with \eqref{finalp}, this last identity implies that $\widehat{q}_{*} = 0$.
\end{remark}

Under some particular choices of the admissible domain and the data of system \eqref{nsstokes0}, the proof of Theorem \ref{epslevel} becomes elementary, as an exact solution is given by the Hagen-Poiseuille flow.

\begin{remark} \label{remexact1}
	Suppose that $\Omega \doteq (0,\ell)\times (-h,h)\subset \R^2$, for some constants $\ell,h>0$. Recalling \eqref{poi1}, the steady-state Hagen-Poiseuille flow associated to $\Omega$, having flow rate $\Phi \geq 0$, is given by
	\begin{equation} \label{exactpoi2dfor}
		\bs{U}_{\! \Phi}(x,y) \doteq \left(\dfrac{3\Phi}{4 h^{3}}(h^2 - y^2)\, , \, 0 \right)  \qquad \text{and} \qquad Q_{\Phi}(x,y) \doteq -\dfrac{3 \nu \Phi}{2 h^{3}} \left( x - \ell \right) \qquad \forall (x,y) \in \overline{\Omega} \, ,
	\end{equation}
	see Figure \ref{poi2d} below for a depiction of the velocity profile in \eqref{exactpoi2dfor}.
	\begin{figure}[H]
		\begin{center}
			\includegraphics[scale=0.52]{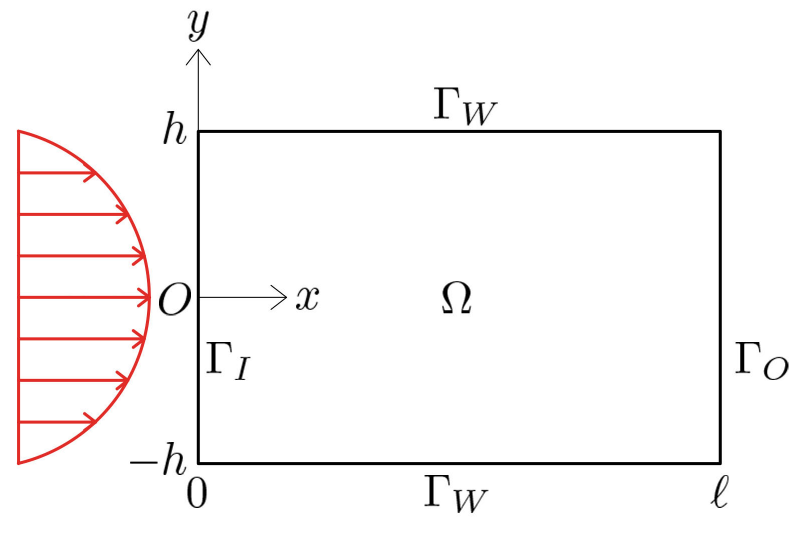}
		\end{center}
		\vspace*{-5mm}
		\caption{Representation of the two-dimensional Hagen-Poiseuille \eqref{exactpoi2dfor}.}\label{poi2d}
	\end{figure}
	\noindent
	A straightforward computation shows that the pair $(\bs{U}_{\! \Phi}, Q_{\Phi}) \in \mathcal{C}^{2}(\overline{\Omega}; \mathbb{R}^{2}) \times \mathcal{C}^{1}(\overline{\Omega}; \mathbb{R})$ solves, in classical sense, the system 
	\begin{equation} \label{exactpoi2d}
		\left\{
		\begin{aligned}
			& -\nu\Delta \bs{U}_{\! \Phi} + (\bs{U}_{\! \Phi} \cdot \nabla)\bs{U}_{\! \Phi}+\nabla Q_{\Phi}=\bs{0} \, , \quad  \nabla\cdot \bs{U}_{\! \Phi}=0 \ \ \mbox{ in } \ \ \Omega \, , \\[5pt]
			& \bs{U}_{\! \Phi}=\bs{g}_{*} \ \ \mbox{ on } \ \ \Gamma_{I} \, , \\[5pt]
			& \bs{U}_{\! \Phi}=\bs{0} \ \ \mbox{ on } \ \ \Gamma_{W} \, , \\[5pt]
			& \nu \dfrac{\partial \bs{U}_{\! \Phi}}{\partial \bs{n}} - Q_{\Phi} \, \bs{n} = \bs{0} \ \ \mbox{ on } \ \ \Gamma_{O} \, ,
		\end{aligned}
		\right.
	\end{equation}
	where $\bs{g}_{*}(x,y) \doteq \bs{U}_{\! \Phi}(x,y)$, for every $(x,y) \in \Gamma_{I}$. Moreover, since
	$$
	\bs{U}_{\! \Phi}(x,y) \cdot \bs{n} = \dfrac{3\Phi}{4 h^{3}}(h^2 - y^2) \geq 0 \qquad \forall (x,y) \in \Gamma_{O} \, ,
	$$
	we deduce that the pair $(\bs{U}_{\! \Phi}, Q_{\Phi})$ not only serves as a reference flow (in the sense of Lemma \ref{lemma0}), it also forms a classical solution to system \eqref{nsstokes0} for the particular choice of data presented in \eqref{exactpoi2d}.
\end{remark}

\begin{remark} \label{remexact2}
	Suppose that $\Omega \doteq \{(x,y,z) \in \mathbb{R}^{3} \ | \ x^2 + y^2 < R^2 \, , \quad 0 < z < \ell \ \}$, for some constants $R,\ell>0$. Recalling \eqref{poi1}, the steady-state Hagen-Poiseuille flow in $\Omega$, with flow rate $\Phi \geq 0$, is given by
	\begin{equation} \label{exactpoi3dfor}
		\bs{U}_{\! \Phi}(x,y,z) \doteq \left(0 \, , 0 \, , \dfrac{2\Phi}{\pi R^{4}}(R^2 - x^{2} - y^{2}) \right)  \quad \text{and} \quad Q_{\Phi}(x,y,z) \doteq -\dfrac{8 \nu \Phi}{\pi R^{4}} \left( z - \ell \right) \quad \forall (x,y,z) \in \overline{\Omega} \, ,
	\end{equation}
	see Figure \ref{poi3d} below for a depiction of the velocity profile in \eqref{exactpoi3dfor}.
	\begin{figure}[H]
		\begin{center}
			\includegraphics[scale=0.57]{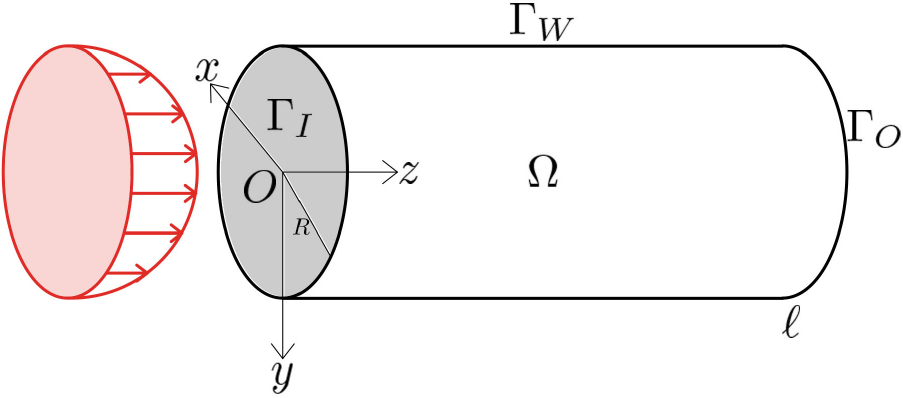}
		\end{center}
		\vspace*{-5mm}
		\caption{Representation of the three-dimensional Hagen-Poiseuille \eqref{exactpoi3dfor}.}\label{poi3d}
	\end{figure}
	\noindent
	A straightforward computation shows that the pair $(\bs{U}_{\! \Phi}, Q_{\Phi}) \in \mathcal{C}^{2}(\overline{\Omega}; \mathbb{R}^{3}) \times \mathcal{C}^{1}(\overline{\Omega}; \mathbb{R})$ solves, in classical sense, the system 
	\begin{equation} \label{exactpoi3d}
		\left\{
		\begin{aligned}
			& -\nu\Delta \bs{U}_{\! \Phi} + (\bs{U}_{\! \Phi} \cdot \nabla)\bs{U}_{\! \Phi}+\nabla Q_{\Phi}=\bs{0} \, , \quad  \nabla\cdot \bs{U}_{\! \Phi}=0 \ \ \mbox{ in } \ \ \Omega \, , \\[5pt]
			& \bs{U}_{\! \Phi}=\bs{g}_{*} \ \ \mbox{ on } \ \ \Gamma_{I} \, , \\[5pt]
			& \bs{U}_{\! \Phi}=\bs{0} \ \ \mbox{ on } \ \ \Gamma_{W} \, , \\[5pt]
			& \nu \dfrac{\partial \bs{U}_{\! \Phi}}{\partial \bs{n}} - Q_{\Phi} \, \bs{n} = \bs{0} \ \ \mbox{ on } \ \ \Gamma_{O} \, ,
		\end{aligned}
		\right.
	\end{equation}
	where $\bs{g}_{*}(x,y,z) \doteq \bs{U}_{\! \Phi}(x,y,z)$, for every $(x,y,z) \in \Gamma_{I}$. Moreover, since
	$$
	\bs{U}_{\! \Phi}(x,y,z) \cdot \bs{n} = \dfrac{2\Phi}{\pi R^{4}}(R^2 - x^{2} - y^{2}) \geq 0 \qquad \forall (x,y,z) \in \Gamma_{O} \, ,
	$$
	we deduce that the pair $(\bs{U}_{\! \Phi}, Q_{\Phi})$ not only serves as a reference flow (in the sense of Lemma \ref{lemma0}), it also forms a classical solution to system \eqref{nsstokes0} for the particular choice of data presented in \eqref{exactpoi3d}.
\end{remark}

\section{Unique solvability of the boundary-value problem} \label{solunica}
The main purpose of this section is to prove that, under a suitable smallness assumption on the data of problem \eqref{nsstokes0}, there exists exactly one weak solution to such system. Following the path adopted in \cite[Section 5]{korobkov2020solvability}, we firstly show that if the data of problem \eqref{nsstokes0} are sufficiently small, then any solution to system \eqref{nsstokes0} can be suitably compared with the reference flow built in Lemma \ref{lemma0}.

\begin{lemma} \label{lemmaunique1}
Let $\Omega \subset \mathbb{R}^{N}$, $N \in \{2,3\}$, be an admissible domain. Given an external force $\bs{f} \in L^2(\Omega;\mathbb{R}^N)$, $\bs{g}_{*} \in W^{3/2,2}_{+}(\Gamma_{I};\mathbb{R}^N)$ and $\sigma_{*} \in W^{1/2,2}(\Gamma_{O};\mathbb{R})$, suppose that $\bs{u}\in \mathcal{V}(\Omega)$ is a weak solution of problem \eqref{nsstokes0}, in the sense of Definition \ref{weaksolution}. There exists $\omega_{*} > 0$, depending exclusively on $\Omega$ and $\nu$, such that, if
\begin{equation} \label{cond0}
\| \bs{g}_{*} \|_{W^{3/2,2}(\Gamma_{I})} < \omega_{*} \, ,
\end{equation}
then there holds the bound
\begin{equation} \label{cond1}
	\| \nabla (\bs{u} - \bs{W}_{\! \! *}) \|_{L^{2}(\Omega)} < C_{*} \left( \| \bs{f} \|_{L^{2}(\Omega)} + \| \sigma_{*} \|_{L^{2}(\Gamma_{O})} + \| \bs{g}_{*} \|_{W^{3/2,2}(\Gamma_{I})} + \| \bs{g}_{*} \|^{2}_{W^{3/2,2}(\Gamma_{I})} \right) \, ,
\end{equation}
for some constant $C_{*} > 0$ that depends only on $\Omega$ and $\nu$.
\end{lemma}
\noindent
\begin{proof}
Denote by $S_{*}>0$ be the \textit{best Sobolev constant} of the embedding $\mathcal{V}_{*}(\Omega) \subset L^{4}(\Omega;\mathbb{R}^N)$, namely
$$
S_{*} \doteq \min_{\bs{v} \in \mathcal{V}_{*}(\Omega) \setminus \{0\}} \ \dfrac{\|\nabla \bs{v}\|^2_{L^2(\Omega)}}{\|\bs{v}\|^{2}_{L^{4}(\Omega)}} \, ,
$$ 
which depends only on $\Omega$ (see \cite{talenti1976best} for further details), so that 
\begin{equation} \label{cond22}
\|\bs{v}\|^{2}_{L^{4}(\Omega)} \leq \dfrac{1}{S_{*}} \|\nabla \bs{v}\|^2_{L^2(\Omega)} \qquad \forall \bs{v} \in \mathcal{V}_{*}(\Omega) \, .
\end{equation}
Letting $M_{*} > 0$ be the constant appearing in the right-hand side of the first inequality in \eqref{fluxestimate}, put
\begin{equation} \label{cond2}
\omega_{*} \doteq \dfrac{\nu S_{*}}{2 M_{*}} \, .
\end{equation}

Define $\widehat{\bs{u}} \doteq \bs{u} - \bs{W}_{\! \! *}$ in $\Omega$, so that $\widehat{\bs{u}} \in \mathcal{V}_{*}(\Omega)$ satisfies the identity
\begin{equation}\label{oseenquasiuni}
	\begin{aligned}
		& \nu \int_{\Omega} \nabla \widehat{\bs{u}} : \nabla \bs{\varphi} + \int_{\Omega} (\widehat{\bs{u}} \cdot \nabla)\widehat{\bs{u}} \cdot \bs{\varphi} + \int_{\Omega} \left[ (\widehat{\bs{u}} \cdot \nabla)\bs{W}_{\! \! *} + (\bs{W}_{\! \! *} \cdot \nabla)\widehat{\bs{u}} \right] \cdot \bs{\varphi} + \dfrac{1}{2} \int_{\Gamma_{O}} [(\widehat{\bs{u}} + \bs{W}_{\! \! *}) \cdot \bs{n}]^{-}(\widehat{\bs{u}} \cdot \bs{\varphi}) \\[6pt]
		& \hspace{-3mm}   = \int_{\Omega} \bs{f} \cdot \bs{\varphi} + \int_{\Gamma_{O}} \sigma_{*} (\bs{\varphi} \cdot \bs{n}) - \nu \int_{\Omega} \nabla \bs{W}_{\! \! *} : \nabla \bs{\varphi} - \int_{\Omega} (\bs{W}_{\! \! *} \cdot \nabla)\bs{W}_{\! \! *} \cdot \bs{\varphi} \qquad \forall \bs{\varphi} \in \mathcal{V}_{*}(\Omega) \, .
	\end{aligned}
\end{equation}
see again \eqref{oseenquasi}. As in \eqref{vlambda1}, an integration by parts yields
\begin{equation} \label{vlambda1uni}
	\begin{aligned}
		\int_{\Omega} (\widehat{\bs{u}} \cdot \nabla)\widehat{\bs{u}} \cdot \widehat{\bs{u}} = \dfrac{1}{2} \int_{\Gamma_{O}} | \widehat{\bs{u}} |^{2}(\widehat{\bs{u}} \cdot \bs{n}) \quad \text{and} \quad \int_{\Omega} (\bs{W}_{\! \! *} \cdot \nabla)\widehat{\bs{u}} \cdot \widehat{\bs{u}} = \dfrac{1}{2} \int_{\Gamma_{O}} | \widehat{\bs{u}} |^{2}(\bs{W}_{\! \! *} \cdot \bs{n}) \, .
	\end{aligned}
\end{equation}
After putting $\bs{\varphi} = \widehat{\bs{u}}$ in \eqref{oseenquasiuni} and enforcing \eqref{vlambda1uni}, the following identity is obtained:
\begin{equation} \label{importantunipre}
\begin{aligned}
	\nu \| \nabla \widehat{\bs{u}} \|^{2}_{L^{2}(\Omega)} & = \int_{\Omega} \bs{f} \cdot \widehat{\bs{u}} + \int_{\Gamma_{O}} \sigma_{*} (\widehat{\bs{u}} \cdot \bs{n}) - \nu \int_{\Omega} \nabla \bs{W}_{\! \! *} : \nabla \widehat{\bs{u}} - \int_{\Omega} (\bs{W}_{\! \! *} \cdot \nabla)\bs{W}_{\! \! *} \cdot \widehat{\bs{u}}  \\[6pt]
	& \hspace{4mm} - \int_{\Omega} (\widehat{\bs{u}} \cdot \nabla)\bs{W}_{\! \! *} \cdot \widehat{\bs{u}} - \dfrac{1}{2} \int_{\Gamma_{O}} \left( (\widehat{\bs{u}} + \bs{W}_{\! \! *}) \cdot \bs{n} + [(\widehat{\bs{u}} + \bs{W}_{\! \! *}) \cdot \bs{n}]^{-} \right) | \widehat{\bs{u}} |^{2} \, .
\end{aligned}
\end{equation}
By successive applications of H\"older, Poincaré, Sobolev (see \eqref{cond22}) and trace inequalities on the right-hand side of  \eqref{importantunipre}, together with \eqref{fluxestimate}, we can further bound in the following way:
\begin{equation} \label{importantuni}
	\begin{aligned}	
		\nu \| \nabla \widehat{\bs{u}} \|^{2}_{L^{2}(\Omega)} & \leq \int_{\Omega} \bs{f} \cdot \widehat{\bs{u}} + \int_{\Gamma_{O}} \sigma_{*} (\widehat{\bs{u}} \cdot \bs{n}) - \nu \int_{\Omega} \nabla \bs{W}_{\! \! *} : \nabla \widehat{\bs{u}} - \int_{\Omega} (\bs{W}_{\! \! *} \cdot \nabla)\bs{W}_{\! \! *} \cdot \widehat{\bs{u}} - \int_{\Omega} (\widehat{\bs{u}} \cdot \nabla)\bs{W}_{\! \! *} \cdot \widehat{\bs{u}} \\[6pt]
		& \leq C \left( \| \bs{f} \|_{L^{2}(\Omega)} + \| \sigma_{*} \|_{L^{2}(\Gamma_{O})} + \| \bs{g}_{*} \|_{W^{3/2,2}(\Gamma_{I})} + \| \bs{g}_{*} \|^{2}_{W^{3/2,2}(\Gamma_{I})} \right) \| \nabla \widehat{\bs{u}} \|_{L^{2}(\Omega)} \\[6pt]
		& \hspace{4mm} + \dfrac{M_{*}}{S_{*}} \| \bs{g}_{*} \|_{W^{3/2,2}(\Gamma_{I})} \| \nabla \widehat{\bs{u}} \|^{2}_{L^{2}(\Omega)} \, ,
	\end{aligned} 
\end{equation}
that is,
\begin{equation} \label{importantuni3}
		\left( \nu - \dfrac{M_{*}}{S_{*}} \| \bs{g}_{*} \|_{W^{3/2,2}(\Gamma_{I})} \right) \| \nabla \widehat{\bs{u}} \|_{L^{2}(\Omega)} \leq C \left( \| \bs{f} \|_{L^{2}(\Omega)} + \| \sigma_{*} \|_{L^{2}(\Gamma_{O})} + \| \bs{g}_{*} \|_{W^{3/2,2}(\Gamma_{I})} + \| \bs{g}_{*} \|^{2}_{W^{3/2,2}(\Gamma_{I})} \right) \, ,
\end{equation}
from where the desired estimate \eqref{cond1} is readily obtained, under the assumption \eqref{cond0} and with the choice of $\omega_{*}>0$ given in \eqref{cond2}.
\end{proof}

We are now in position to prove the main result of this section.

\begin{theorem} \label{theounique1}
	Let $\Omega \subset \mathbb{R}^{N}$, $N \in \{2,3\}$, be an admissible domain. Let $\bs{f} \in L^2(\Omega;\mathbb{R}^N)$ be an external force, $\bs{g}_{*} \in W^{3/2,2}_{+}(\Gamma_{I};\mathbb{R}^N)$ and $\sigma_{*} \in W^{1/2,2}(\Gamma_{O};\mathbb{R})$. There exists $\varrho_{*} > 0$, depending exclusively on $\Omega$ and $\nu$, such that, if
	\begin{equation} \label{cond0theo}
	\| \bs{f} \|_{L^{2}(\Omega)} +  \| \bs{g}_{*} \|_{W^{3/2,2}(\Gamma_{I})} + \| \sigma_{*} \|_{L^{2}(\Gamma_{O})} < \varrho_{*} \, ,
	\end{equation}
	then the problem \eqref{nsstokes0} admits a unique weak solution in the sense of Definition \ref{weaksolution}.
\end{theorem}
\noindent
\begin{proof}	
Under the assumptions of the statement, the existence of (at least) one weak solution $\bs{u}\in \mathcal{V}(\Omega)$ to problem \eqref{nsstokes0} is ensured by Theorem \ref{epslevel}. Now, suppose there exists another weak solution $\bs{v}\in \mathcal{V}(\Omega)$ to problem \eqref{nsstokes0}. Define $\widehat{\bs{u}} \doteq \bs{u} - \bs{W}_{\! \! *}$ and $\widehat{\bs{v}} \doteq \bs{v} - \bs{W}_{\! \! *}$ in $\Omega$, so that $\widehat{\bs{u}}, \widehat{\bs{v}} \in \mathcal{V}_{*}(\Omega)$ satisfy the identities
\begin{equation}\label{oseenquasiuni1}
	\begin{aligned}
		& \nu \int_{\Omega} \nabla \widehat{\bs{u}} : \nabla \bs{\varphi} + \int_{\Omega} (\widehat{\bs{u}} \cdot \nabla)\widehat{\bs{u}} \cdot \bs{\varphi} + \int_{\Omega} \left[ (\widehat{\bs{u}} \cdot \nabla)\bs{W}_{\! \! *} + (\bs{W}_{\! \! *} \cdot \nabla)\widehat{\bs{u}} \right] \cdot \bs{\varphi} + \dfrac{1}{2} \int_{\Gamma_{O}} [(\widehat{\bs{u}} + \bs{W}_{\! \! *}) \cdot \bs{n}]^{-}(\widehat{\bs{u}} \cdot \bs{\varphi}) \\[6pt]
		& \hspace{-3mm}   = \int_{\Omega} \bs{f} \cdot \bs{\varphi} + \int_{\Gamma_{O}} \sigma_{*} (\bs{\varphi} \cdot \bs{n}) - \nu \int_{\Omega} \nabla \bs{W}_{\! \! *} : \nabla \bs{\varphi} - \int_{\Omega} (\bs{W}_{\! \! *} \cdot \nabla)\bs{W}_{\! \! *} \cdot \bs{\varphi} \qquad \forall \bs{\varphi} \in \mathcal{V}_{*}(\Omega) 
	\end{aligned}
\end{equation}
and 
\begin{equation}\label{oseenquasiuni2}
	\begin{aligned}
		& \nu \int_{\Omega} \nabla \widehat{\bs{v}} : \nabla \bs{\varphi} + \int_{\Omega} (\widehat{\bs{v}} \cdot \nabla)\widehat{\bs{v}} \cdot \bs{\varphi} + \int_{\Omega} \left[ (\widehat{\bs{v}} \cdot \nabla)\bs{W}_{\! \! *} + (\bs{W}_{\! \! *} \cdot \nabla)\widehat{\bs{v}} \right] \cdot \bs{\varphi} + \dfrac{1}{2} \int_{\Gamma_{O}} [(\widehat{\bs{v}} + \bs{W}_{\! \! *}) \cdot \bs{n}]^{-}(\widehat{\bs{v}} \cdot \bs{\varphi}) \\[6pt]
		& \hspace{-3mm}   = \int_{\Omega} \bs{f} \cdot \bs{\varphi} + \int_{\Gamma_{O}} \sigma_{*} (\bs{\varphi} \cdot \bs{n}) - \nu \int_{\Omega} \nabla \bs{W}_{\! \! *} : \nabla \bs{\varphi} - \int_{\Omega} (\bs{W}_{\! \! *} \cdot \nabla)\bs{W}_{\! \! *} \cdot \bs{\varphi} \qquad \forall \bs{\varphi} \in \mathcal{V}_{*}(\Omega) \, .
	\end{aligned}
\end{equation}
see \eqref{oseenquasiuni}. Define now $\bs{z} \doteq \widehat{\bs{v}} - \widehat{\bs{u}} =  \bs{v} - \bs{u}$ in $\Omega$, so that $\bs{z} \in \mathcal{V}_{*}(\Omega)$, and by taking the difference between the integral identities \eqref{oseenquasiuni1}-\eqref{oseenquasiuni2}, we observe that
\begin{equation}\label{oseenquasiuni3}
	\begin{aligned}
		& \nu \int_{\Omega} \nabla \bs{z} : \nabla \bs{\varphi} + \int_{\Omega} \left[ (\widehat{\bs{v}} \cdot \nabla)\widehat{\bs{v}} - (\widehat{\bs{u}} \cdot \nabla)\widehat{\bs{u}} \right] \cdot \bs{\varphi} + \int_{\Omega} \left[ (\bs{z} \cdot \nabla)\bs{W}_{\! \! *} + (\bs{W}_{\! \! *} \cdot \nabla)\bs{z} \right] \cdot \bs{\varphi}  \\[6pt]
		& \hspace{-3mm}  = \dfrac{1}{2} \int_{\Gamma_{O}} [(\widehat{\bs{u}} + \bs{W}_{\! \! *}) \cdot \bs{n}]^{-}(\widehat{\bs{u}} \cdot \bs{\varphi}) - \dfrac{1}{2} \int_{\Gamma_{O}} [(\widehat{\bs{v}} + \bs{W}_{\! \! *}) \cdot \bs{n}]^{-}(\widehat{\bs{v}} \cdot \bs{\varphi}) \qquad \forall \bs{\varphi} \in \mathcal{V}_{*}(\Omega) \, .
	\end{aligned}
\end{equation}
We take $\bs{\varphi} = \bs{z}$ in \eqref{oseenquasiuni3}, and recalling the first identity in \eqref{vlambda1uni}, we obtain
\begin{equation} \label{importantunitheo}
	\begin{aligned}	
		\nu \| \nabla \bs{z} \|^{2}_{L^{2}(\Omega)} & = - \int_{\Omega} \left[ (\bs{z} \cdot \nabla)\widehat{\bs{u}} + (\widehat{\bs{u}} \cdot \nabla)\bs{z} \right] \cdot \bs{z} - \int_{\Omega} \left[ (\bs{z} \cdot \nabla)\bs{W}_{\! \! *} + (\bs{W}_{\! \! *} \cdot \nabla)\bs{z} \right] \cdot \bs{z}  \\[6pt]
		& \hspace{4mm} + \dfrac{1}{2} \int_{\Gamma_{O}} [\bs{u} \cdot \bs{n}]^{-}(\widehat{\bs{u}} \cdot \bs{z}) - \dfrac{1}{2} \int_{\Gamma_{O}} [\bs{v} \cdot \bs{n}]^{-}(\widehat{\bs{v}} \cdot \bs{z}) - \dfrac{1}{2} \int_{\Gamma_{O}} | \bs{z} |^{2}(\bs{z} \cdot \bs{n})  \, .
	\end{aligned} 
\end{equation}
Observe that, almost everywhere on $\Gamma_{O}$, there holds the following:
$$
\begin{aligned}	
	& [\bs{u} \cdot \bs{n}]^{-}(\widehat{\bs{u}} \cdot \bs{z}) - [\bs{v} \cdot \bs{n}]^{-}(\widehat{\bs{v}} \cdot \bs{z}) - (\bs{z} \cdot \bs{n}) | \bs{z} |^{2} \\[6pt]
	& \hspace{-5mm} = \left( [\bs{u} \cdot \bs{n}]^{-} - [\bs{v} \cdot \bs{n}]^{-}  \right) (\widehat{\bs{u}} \cdot \bs{z}) + (\bs{u} \cdot \bs{n}) | \bs{z} |^{2} - \left( [\bs{v} \cdot \bs{n}]^{-} + (\bs{v} \cdot \bs{n}) \right) | \bs{z} |^{2} \\[6pt]
	& \hspace{-5mm} \leq \left( [\bs{u} \cdot \bs{n}]^{-} - [\bs{v} \cdot \bs{n}]^{-}  \right) (\widehat{\bs{u}} \cdot \bs{z}) + (\bs{u} \cdot \bs{n}) | \bs{z} |^{2} \, ,
\end{aligned} 
$$
which, once inserted into the right-hand side of \eqref{importantunitheo}, entails
\begin{equation} \label{importantunitheo2}
	\begin{aligned}	
		\nu \| \nabla \bs{z} \|^{2}_{L^{2}(\Omega)} & \leq - \int_{\Omega} \left[ (\bs{z} \cdot \nabla)\widehat{\bs{u}} + (\widehat{\bs{u}} \cdot \nabla)\bs{z} \right] \cdot \bs{z} - \int_{\Omega} \left[ (\bs{z} \cdot \nabla)\bs{W}_{\! \! *} + (\bs{W}_{\! \! *} \cdot \nabla)\bs{z} \right] \cdot \bs{z}  \\[6pt]
		& \hspace{4mm} + \dfrac{1}{2} \int_{\Gamma_{O}} \left( [\bs{u} \cdot \bs{n}]^{-} - [\bs{v} \cdot \bs{n}]^{-}  \right) (\widehat{\bs{u}} \cdot \bs{z}) + \dfrac{1}{2} \int_{\Gamma_{O}} (\bs{u} \cdot \bs{n}) | \bs{z} |^{2}  \, .
	\end{aligned} 
\end{equation}
By the H\"older, Poincaré and Sobolev inequalities, together with \eqref{fluxestimate}, we obtain
\begin{equation} \label{importantunitheo3}
	\begin{aligned}	
& \left| \int_{\Omega} \left[ (\bs{z} \cdot \nabla)\widehat{\bs{u}} + (\widehat{\bs{u}} \cdot \nabla)\bs{z} \right] \cdot \bs{z} + \int_{\Omega} \left[ (\bs{z} \cdot \nabla)\bs{W}_{\! \! *} + (\bs{W}_{\! \! *} \cdot \nabla)\bs{z} \right] \cdot \bs{z} \right|  \\[6pt]
		& \hspace{-4mm} \leq C \left( \| \nabla \widehat{\bs{u}} \|_{L^{2}(\Omega)} + \| \bs{g}_{*} \|_{W^{3/2,2}(\Gamma_{I})} \right) \| \nabla \bs{z} \|^{2}_{L^{2}(\Omega)}  \, .
	\end{aligned} 
\end{equation}
On the other hand, since 
$$
\left| [r]^{-} - [s]^{-} \right| \leq |r-s| \qquad \forall r,s \in \mathbb{R} \, , 
$$
we may employ the trace inequality (see Remark \ref{rem1}) to estimate the boundary integrals appearing on the right-hand side of \eqref{importantunitheo2} in the following way:
\begin{equation} \label{importantunitheo4}
	\begin{aligned}	
		& \left| \int_{\Gamma_{O}} \left( [\bs{u} \cdot \bs{n}]^{-} - [\bs{v} \cdot \bs{n}]^{-}  \right) (\widehat{\bs{u}} \cdot \bs{z}) + \int_{\Gamma_{O}} (\bs{u} \cdot \bs{n}) | \bs{z} |^{2} \right| \leq \int_{\Gamma_{O}} |\bs{z} \cdot \bs{n}| |\widehat{\bs{u}} \cdot \bs{z}| + \int_{\Gamma_{O}} |\bs{u} \cdot \bs{n}| | \bs{z} |^{2} \\[6pt]	
		& \hspace{-4mm} \leq \left( \| \bs{u} \|_{L^{2}(\Gamma_{O})} + \| \widehat{\bs{u}} \|_{L^{2}(\Gamma_{O})} \right) \| \bs{z} \|^{2}_{L^{4}(\Gamma_{O})} \leq C \left( \| \nabla \widehat{\bs{u}} \|_{L^{2}(\Omega)} + \| \bs{g}_{*} \|_{W^{3/2,2}(\Gamma_{I})} \right) \| \nabla \bs{z} \|^{2}_{L^{2}(\Omega)}  \, .
	\end{aligned} 
\end{equation}
Upon insertion of \eqref{importantunitheo3}-\eqref{importantunitheo4} into the right-hand side of \eqref{importantunitheo2}, we deduce that
\begin{equation} \label{importantunitheo5}
\| \nabla \bs{z} \|^{2}_{L^{2}(\Omega)} \leq C \left( \| \nabla \widehat{\bs{u}} \|_{L^{2}(\Omega)} + \| \bs{g}_{*} \|_{W^{3/2,2}(\Gamma_{I})} \right) \| \nabla \bs{z} \|^{2}_{L^{2}(\Omega)} \, .
\end{equation}

Now, let $\omega_{*} > 0$ be as in the statement of Lemma \ref{lemmaunique1}. For any $\varrho_{*} \in (0,\omega_{*})$, under the assumption \eqref{cond0theo} we can certainly claim that the condition \eqref{cond0} is fulfilled, implying that \eqref{cond1} holds. Inserted into the right-hand side of \eqref{importantunitheo5}, this means that
\begin{equation} \label{importantunitheo6}
	\| \nabla \bs{z} \|^{2}_{L^{2}(\Omega)} \leq C \left( \| \bs{f} \|_{L^{2}(\Omega)} + \| \sigma_{*} \|_{L^{2}(\Gamma_{O})} + \| \bs{g}_{*} \|_{W^{3/2,2}(\Gamma_{I})} + \| \bs{g}_{*} \|^{2}_{W^{3/2,2}(\Gamma_{I})} \right) \| \nabla \bs{z} \|^{2}_{L^{2}(\Omega)} \, .
\end{equation}
Thus, denoting by $M_{*} > 0$ the constant entering the right-hand side of \eqref{importantunitheo6}, and choosing $\varrho_{*} \in (0,\omega_{*})$ in such a way that $M_{*} (\varrho_{*} + \varrho^{2}_{*}) < 1$, we reach the desired conclusion under the assumption \eqref{cond0theo}.
\end{proof}

\newpage
\begin{remark}
	All the results presented in this manuscript can be easily adapted and extended to the case when the admissible domain $\Omega \subset \mathbb{R}^N$ has several rectangular/cylindrical inlets and outlets. Examples of such junctions of pipes are depicted in Figures \ref{polpo2d}-\ref{polpo3d}.
	\begin{figure}[H]
		\begin{center}
			\includegraphics[scale=0.52]{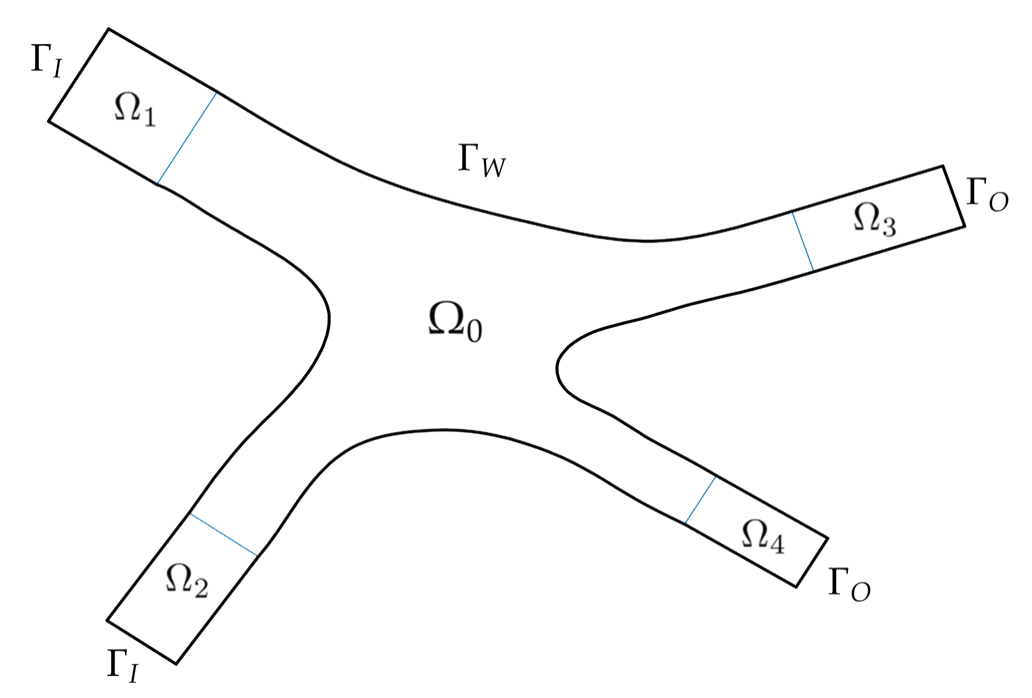}
		\end{center}
		\vspace*{-6mm}
		\caption{Representation of an admissible domain $\Omega \subset \mathbb{R}^{2}$ with several inlets and outlets.}\label{polpo2d}
	\end{figure}
	
	\begin{figure}[H]
		\begin{center}
			\includegraphics[scale=0.5]{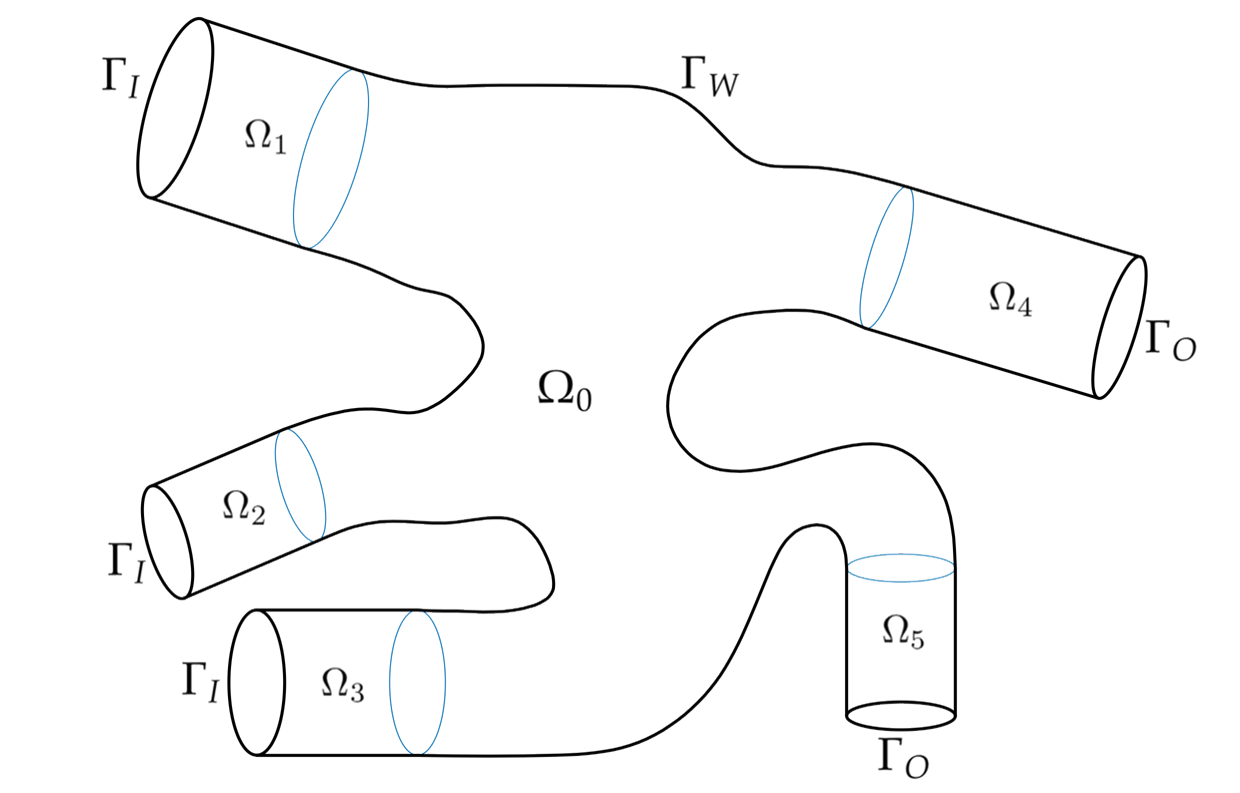}
		\end{center}
		\vspace*{-6mm}
		\caption{Representation of an admissible domain $\Omega \subset \mathbb{R}^{3}$ with several inlets and outlets.}\label{polpo3d}
	\end{figure}
\end{remark}

\bigskip

\noindent
{\bf Acknowledgements.} 
The research of Alessio Falocchi is supported by the grant \textit{Dipartimento di Eccellenza 2023-2027},
issued by the Ministry of University and Research (Italy) and is a part of the INdAM-GNAMPA project entitled ``EDP e Applicazioni: Dinamica dei Fluidi e Teoria Spettrale'' (15/01/2026-31/12/2026). Ana L. Silvestre acknowledges the financial support of \textit{Funda\c{c}\~ao para a Ci\^encia e a Tecnologia}  (FCT, Portuguese Agency for Scientific Research), through the project UIDB/04621/2025 of CEMAT/IST-ID (https://doi.org/10.54499/UID/04621/2025). The research of Gianmarco Sperone is supported by the \textit{Chilean National Agency for Research and Development} (ANID) through the \textit{Fondecyt Iniciación} grant 11250322. The Authors warmly thank two anonymous Referees for their careful proofreading and for several suggestions that led to an improvement of the present work.
\par\smallskip
\noindent
{\bf Data availability statement.} Data sharing not applicable to this article as no datasets were generated or analyzed during the current study.
\par\smallskip
\noindent
{\bf Conflict of interest statement}.  The Authors declare that they have no conflict of interest.

\phantomsection
\addcontentsline{toc}{section}{References}
\bibliographystyle{abbrv}
\bibliography{references}

\begin{thebibliography}{10}

\bibitem{amick1977steady}
C.~J. Amick.
\newblock Steady solutions of the {N}avier-{S}tokes equations in unbounded
  channels and pipes.
\newblock {\em Annali della Scuola Normale Superiore di Pisa - Classe di
  Scienze}, 4(3):473--513, 1977.

\bibitem{amick1984existence}
C.~J. Amick.
\newblock Existence of solutions to the nonhomogeneous steady {N}avier-{S}tokes
  equations.
\newblock {\em Indiana University Mathematics Journal}, 33(6):817--830, 1984.

\bibitem{Baranovskii2020OptimalBC}
E.~S. Baranovskii.
\newblock Optimal boundary control of nonlinear-viscous fluid flows.
\newblock {\em Sbornik: Mathematics}, 211:505 -- 520, 2020.

\bibitem{benevs2014note}
M.~Bene{\v{s}}.
\newblock A note on regularity and uniqueness of natural convection with
  effects of viscous dissipation in 3{D} open channels.
\newblock {\em Zeitschrift f{\"u}r angewandte Mathematik und Physik},
  65:961--975, 2014.

\bibitem{benevs2016solutions}
M.~Bene{\v{s}} and P.~Ku{\v{c}}era.
\newblock Solutions to the {N}avier--{S}tokes equations with mixed boundary
  conditions in two-dimensional bounded domains.
\newblock {\em Mathematische Nachrichten}, 289(2-3):194--212, 2016.

\bibitem{BCBBG2018}
C.~Bertoglio, A.~Caiazzo, Y.~Bazilevs, M.~Braack, M.~Esmaily, V.~Gravemeier,
  A.~L.~Marsden, O.~Pironneau, I.~E.~Vignon-Clementel, and W.~A.~Wall.
\newblock Benchmark problems for numerical treatment of backflow at open
  boundaries.
\newblock {\em International Journal for Numerical Methods in Biomedical
  Engineering}, 34(2):e2918, 2018.

\bibitem{blazy2007artificial}
S.~Blazy, S.~A. Nazarov, and M.~Specovius-Neugebauer.
\newblock Artificial boundary conditions of pressure type for viscous flows in
  a system of pipes.
\newblock {\em Journal of Mathematical Fluid Mechanics}, 9(1):1--33, 2007.

\bibitem{braack2014directional}
M.~Braack and P.~B. Mucha.
\newblock Directional do-nothing condition for the {N}avier-{S}tokes equations.
\newblock {\em Journal of Computational Mathematics}, pages 507--521, 2014.

\bibitem{bradley1987petroleum}
H.~B. Bradley.
\newblock {\em Petroleum {E}ngineering {H}andbook}.
\newblock Society of Petroleum Engineers, Richardson, Texas, 1987.

\bibitem{bruneau2000boundary}
C.-H. Bruneau.
\newblock Boundary conditions on artificial frontiers for incompressible and
  compressible {N}avier-{S}tokes equations.
\newblock {\em ESAIM: Mathematical Modelling and Numerical Analysis},
  34(2):303--314, 2000.

\bibitem{BruneauFabrie1994}
C.-H. Bruneau and P.~Fabrie.
\newblock Effective downstream boundary conditions for incompressible
  {N}avier-{S}tokes equations.
\newblock {\em International Journal for Numerical Methods in Fluids},
  19(8):693--705, 1994.

\bibitem{bruneau1996new}
C.-H. Bruneau and P.~Fabrie.
\newblock New efficient boundary conditions for incompressible
  {N}avier-{S}tokes equations: a well-posedness result.
\newblock {\em ESAIM: Mathematical Modelling and Numerical Analysis},
  30(7):815--840, 1996.

\bibitem{cattabriga1961problema}
L.~Cattabriga.
\newblock Su un problema al contorno relativo al sistema di equazioni di
  {S}tokes.
\newblock {\em Rendiconti del Seminario Matematico della Università di
  Padova}, 31:308--340, 1961.

\bibitem{concaenglish}
C.~Conca, F.~Murat, and O.~Pironneau.
\newblock The {S}tokes and {N}avier-{S}tokes equations with boundary conditions
  involving the pressure.
\newblock {\em Japanese Journal of Mathematics}, 20(2):279--318, 1994.

\bibitem{dautray1999mathematical}
R.~Dautray and J.-L. Lions.
\newblock {\em Mathematical {A}nalysis and {N}umerical {M}ethods for {S}cience
  and {T}echnology - {V}olume 2: {F}unctional and {V}ariational {M}ethods}.
\newblock Springer Science \& Business Media, 1999.

\bibitem{DONG2015300}
S.~Dong.
\newblock A convective-like energy-stable open boundary condition for
  simulations of incompressible flows.
\newblock {\em Journal of Computational Physics}, 302:300--328, 2015.

\bibitem{FMN2007}
L.~Formaggia, A.~Moura, and F.~Nobile.
\newblock On the stability of the coupling of 3{D} and 1{D} fluid-structure
  interaction models for blood flow simulations.
\newblock {\em ESAIM: Mathematical Modelling and Numerical Analysis},
  41(4):743--769, 2007.

\bibitem{fursikov2009optimal}
A.~Fursikov and R.~Rannacher.
\newblock Optimal {N}eumann control for the two-dimensional steady-state
  {N}avier-{S}tokes equations.
\newblock In {\em New Directions in Mathematical Fluid Mechanics}, pages
  193--221. Springer, 2009.

\bibitem{galdi2011introduction}
G.~P. Galdi.
\newblock {\em An {I}ntroduction to the {M}athematical {T}heory of the
  {N}avier-{S}tokes {E}quations: {S}teady-{S}tate {P}roblems}.
\newblock Springer Science \& Business Media, 2011.

\bibitem{galdi2008hemodynamical}
G.~P. Galdi, A.~M. Robertson, R.~Rannacher, and S.~Turek.
\newblock {\em Hemodynamical {F}lows: {M}odeling, {A}nalysis and {S}imulation
  ({O}berwolfach {S}eminars)}.
\newblock Springer Science \& Business Media, 2008.

\bibitem{gazzola2025steady}
F.~Gazzola, M.~V. Korobkov, X.~Ren, and G.~Sperone.
\newblock The steady {N}avier-{S}tokes equations in a system of unbounded
  channels with sources and sinks.
\newblock {\em arXiv preprint arXiv:2505.14642}, 2025.

\bibitem{girault2012finite}
V.~Girault and P.-A. Raviart.
\newblock {\em Finite {E}lement {M}ethods for {N}avier-{S}tokes {E}quations:
  {T}heory and {A}lgorithms}, volume~5.
\newblock Springer Science \& Business Media, 2012.

\bibitem{gresho1991some}
P.~M. Gresho.
\newblock Some current {CFD} issues relevant to the incompressible
  {N}avier-{S}tokes equations.
\newblock {\em Computer Methods in Applied Mechanics and Engineering},
  87(2-3):201--252, 1991.

\bibitem{grisvard2011elliptic}
P.~Grisvard.
\newblock {\em Elliptic {P}roblems in {N}onsmooth {D}omains}.
\newblock Society for Industrial and Applied Mathematics, Philadelphia, 2011.

\bibitem{heywood1996artificial}
J.~G. Heywood, R.~Rannacher, and S.~Turek.
\newblock Artificial boundaries and flux and pressure conditions for the
  incompressible {N}avier--{S}tokes equations.
\newblock {\em International Journal for Numerical Methods in Fluids},
  22(5):325--352, 1996.

\bibitem{john2002higher}
V.~John.
\newblock Higher order finite element methods and multigrid solvers in a
  benchmark problem for the 3{D} {N}avier--{S}tokes equations.
\newblock {\em International Journal for Numerical Methods in Fluids},
  40(6):775--798, 2002.

\bibitem{kapitanskii1983spaces}
L.~V. Kapitanskii and K.~Pileckas.
\newblock Spaces of solenoidal vector fields in boundary value problems for the
  {N}avier-{S}tokes equations in regions with noncompact boundaries.
\newblock {\em Matematicheskii Institut imeni Steklova Trudy}, 159:5--36, 1983.

\bibitem{korobkov2011bernoulli}
M.~V. Korobkov.
\newblock Bernoulli law under minimal smoothness assumptions.
\newblock {\em Doklady Mathematics}, 83(1):107--110, 2011.

\bibitem{korobkov2015solution}
M.~V. Korobkov, K.~Pileckas, and R.~Russo.
\newblock Solution of {L}eray's problem for stationary {N}avier-{S}tokes
  equations in plane and axially symmetric spatial domains.
\newblock {\em Annals of Mathematics}, 181(2):769--807, 2015.

\bibitem{korobkov2020solvability}
M.~V. Korobkov, K.~Pileckas, and R.~Russo.
\newblock Solvability in a finite pipe of steady-state {N}avier--{S}tokes
  equations with boundary conditions involving {B}ernoulli pressure.
\newblock {\em Calculus of Variations and Partial Differential Equations},
  59(1):1--22, 2020.

\bibitem{korobkov2024steady}
M.~V. Korobkov, K.~Pileckas, and R.~Russo.
\newblock {\em The {S}teady {N}avier-{S}tokes {S}ystem: {B}asics of the
  {T}heory and the {L}eray {P}roblem}.
\newblock Springer, 2024.

\bibitem{kravcmar2001weak}
S.~Kra{\v{c}}mar and J.~Neustupa.
\newblock A weak solvability of a steady variational inequality of the
  {N}avier--{S}tokes type with mixed boundary conditions.
\newblock {\em Nonlinear Analysis: Theory, Methods \& Applications},
  47(6):4169--4180, 2001.

\bibitem{kravcmar2018modeling}
S.~Kra{\v{c}}mar and J.~Neustupa.
\newblock Modeling of the unsteady flow through a channel with an artificial
  outflow condition by the {N}avier--{S}tokes variational inequality.
\newblock {\em Mathematische Nachrichten}, 291(11-12):1801--1814, 2018.

\bibitem{ladyzhenskaya1969mathematical}
O.~A. Ladyzhenskaya.
\newblock {\em The {M}athematical {T}heory of {V}iscous {I}ncompressible
  {F}low}, volume~76.
\newblock Gordon and Breach, New York, 1969.

\bibitem{ladyzhenskaia1979determination}
O.~A. Ladyzhenskaya and V.~A. Solonnikov.
\newblock Determination of the solutions of {N}avier-{S}tokes stationary
  boundary-value problems with infinite dissipation in unbounded regions.
\newblock {\em Doklady Akademii Nauk USSR}, 249(4):828--831, 1979.

\bibitem{landau}
L.~Landau and E.~Lifshitz.
\newblock {\em Theoretical {P}hysics: {F}luid {M}echanics}, volume~6.
\newblock Pergamon Press, 1987.

\bibitem{lanzendorfer2020multiple}
M.~Lanzend{\"o}rfer and J.~Hron.
\newblock On multiple solutions to the steady flow of incompressible fluids
  subject to do-nothing or constant traction boundary conditions on artificial
  boundaries.
\newblock {\em Journal of Mathematical Fluid Mechanics}, 22(1):1--18, 2020.

\bibitem{leray1933etude}
J.~Leray.
\newblock {\'E}tude de diverses {\'e}quations int{\'e}grales non lin{\'e}aires
  et de quelques probl{\`e}mes que pose l'hydrodynamique.
\newblock {\em Journal de Math{\'e}matiques Pures et Appliqu{\'e}es}, 12:1--82,
  1933.

\bibitem{mazya2010elliptic}
V.~G. Maz'ya and J.~Rossmann.
\newblock {\em Elliptic {E}quations in {P}olyhedral {D}omains}.
\newblock American Mathematical Society, 2010.

\bibitem{nazarov2008artificial}
S.~A. Nazarov and M.~Specovius-Neugebauer.
\newblock Artificial boundary conditions for the {S}tokes and {N}avier-{S}tokes
  equations in domains that are layer-like at infinity.
\newblock {\em Zeitschrift f{\"u}r Analysis und ihre Anwendungen},
  27(2):125--155, 2008.

\bibitem{necas2011direct}
J.~Ne{\v{c}}as.
\newblock {\em Direct {M}ethods in the {T}heory of {E}lliptic {E}quations}.
\newblock Springer Science \& Business Media, 2011.

\bibitem{neustupa2022maximum}
T.~Neustupa.
\newblock The maximum regularity property of the steady {S}tokes problem
  associated with a flow through a profile cascade in ${L}^{r}$-framework.
\newblock {\em Applications of Mathematics}, 68(2):171--190, 2022.

\bibitem{neustupa2023existence}
T.~Neustupa.
\newblock Existence of a steady flow through a rotating radial turbine with an
  arbitrarily large inflow and an artificial boundary condition on the outflow.
\newblock {\em Zeitschrift f{\"u}r Angewandte Mathematik und Mechanik},
  103(10):e202200439, 2023.

\bibitem{nogueira2025regularized}
P.~Nogueira and A.~L. Silvestre.
\newblock Regularized directional do-nothing boundary conditions for the
  {N}avier-{S}tokes equations: Analytical and numerical study.
\newblock {\em Applied Mathematics and Computation}, 499:129398, 2025.

\bibitem{nogueira2025steady}
P.~Nogueira, A.~L. Silvestre, and J.~Tiago.
\newblock Steady {N}avier--{S}tokes equations with regularized directional
  do-nothing boundary condition: Optimal boundary control for a velocity
  tracking problem.
\newblock {\em Applied Mathematics \& Optimization}, 91(1):24, 2025.

\bibitem{ott2013mixed}
K.~A. Ott and R.~M. Brown.
\newblock The mixed problem for the {L}aplacian in {L}ipschitz domains.
\newblock {\em Potential Analysis}, 38(4):1333--1364, 2013.

\bibitem{pileckas1983spaces}
K.~Pileckas.
\newblock Spaces of solenoidal vectors.
\newblock {\em Trudy Matematicheskogo Instituta imeni V. A. Steklova},
  159:137--149, 1983.

\bibitem{pileckas2007navier}
K.~Pileckas.
\newblock {N}avier--{S}tokes system in domains with cylindrical outlets to
  infinity. {L}eray’s problem.
\newblock In {\em Handbook of Mathematical Fluid Dynamics}, volume~4, pages
  445--647. Elsevier Amsterdam, 2007.

\bibitem{rannacher2012short}
R.~Rannacher.
\newblock A short course on numerical simulation of viscous flow:
  discretization, optimization and stability analysis.
\newblock {\em Discrete \& Continuous Dynamical Systems-S}, 5(6):1147, 2012.

\bibitem{sperone2021steady}
G.~Sperone.
\newblock On the steady motion of {N}avier--{S}tokes flows past a fixed
  obstacle in a three-dimensional channel under mixed boundary conditions.
\newblock {\em Annali di Matematica Pura ed Applicata (1923-)},
  200(5):1961--1985, 2021.

\bibitem{talenti1976best}
G.~Talenti.
\newblock Best constant in {S}obolev inequality.
\newblock {\em Annali di Matematica Pura ed Applicata}, 110:353--372, 1976.

\bibitem{von2004aerodynamics}
T.~von K{\'a}rm{\'a}n.
\newblock {\em Aerodynamics: {S}elected {T}opics in the {L}ight of their
  {H}istorical {D}evelopment}.
\newblock Courier Corporation, 2004.

\bibitem{yosida2012functional}
K.~Yosida.
\newblock {\em Functional {A}nalysis}.
\newblock Springer Science \& Business Media, 1995.

\bibitem{zeidler2013nonlinear}
E.~Zeidler.
\newblock {\em Nonlinear {F}unctional {A}nalysis and its {A}pplications {I}:
  {F}ixed-{P}oint Theorems}.
\newblock Springer Science \& Business Media, 2013.

\end{thebibliography}
\vspace{5mm}

	\noindent
	\hspace{0.1mm}
	\begin{minipage}{140mm}
		\textbf{Alessio Falocchi}\\
		Dipartimento di Matematica\\
		Dipartimento di Eccellenza MUR 2023-2027\\
		Politecnico di Milano\\
		Piazza Leonardo da Vinci 32\\
		20133 Milan - Italy\\
		E-mail: alessio.falocchi@polimi.it
		\vspace{0.5cm}	
	\end{minipage}
	\newline
	\vspace{0.5cm}
	\noindent
	\begin{minipage}{100mm}
		\textbf{Ana Leonor Silvestre}\\
		Centro de Matemática Computacional e Estocástica\\
		Instituto Superior Técnico, Universidade de Lisboa\\
		Avenida Rovisco Pais 1\\
		1049-001 Lisbon - Portugal\\
		and\\
		Departamento de Matemática do Instituto Superior Técnico\\
		Universidade de Lisboa\\
		Avenida Rovisco Pais 1\\
		1049-001 Lisbon - Portugal\\
		E-mail: ana.silvestre@math.tecnico.ulisboa.pt
	\end{minipage}
	\newline
	\vspace{0.5cm}
	\begin{minipage}{100mm}
		\textbf{Gianmarco Sperone}\\
		Facultad de Matemáticas\\
		Pontificia Universidad Católica de Chile\\
		Avenida Vicuña Mackenna 4860\\
		7820436 Santiago - Chile\\
		E-mail: gianmarco.sperone@uc.cl
	\end{minipage}

\end{document}